\documentclass[11pt, oneside]{article}
\usepackage{amsfonts}
\usepackage{mathrsfs}
\usepackage{color}
\usepackage{latexsym}
\usepackage{amssymb}
\usepackage{amsmath}
\usepackage{enumerate}
\usepackage{amsthm}
\usepackage{indentfirst}
\usepackage{mathtools}
\usepackage{tikz}
\usepackage{psfrag}
\usepackage{graphicx}
\usepackage{bm}
\usepackage[rflt]{floatflt}
\usepackage{float}
\usepackage{authblk}
\usepackage{ulem}
\usepackage[square, comma, sort&compress, numbers]{natbib}
\usepackage{verbatim}
\usepackage{amsthm}
\usepackage{esint}

\numberwithin{equation}{section}
\allowdisplaybreaks

\newenvironment{proof2.1}{\medskip\noindent{\bf Proof of the Theorem 2.1:}\enspace}{\hfill \qed \newline \medskip}
\newenvironment{proof2.2}{\medskip\noindent{\bf Proof of the Theorem 2.2:}\enspace}{\hfill \qed \newline \medskip}

\newtheorem{theorem}{\color{black}\indent Theorem}[section]
\newtheorem{lemma}{\color{black}\indent Lemma}[section]

\newtheorem{definition}{\color{black}\indent Definition}[section]
\newtheorem{remark}{\color{black}\indent Remark}[section]

\usepackage{amssymb,amsmath}
\begin{document}
\title{Normalized solutions to  a class of Kirchhoff type equations with a logarithmic perturbation}
\author{Qi Li$^{1}$ \qquad Wenshu Zhou$^{1}$ \qquad
Yuzhu Han$^{2,\dag}$}

\affil{{\it\small $^1$School of Mathematical Sciences, Dalian Minzu University, Dalian 116600, P.R. China}}
 \affil{{\it\small $^2$ School of Mathematics, Jilin University, Changchun 130012, P.R. China}}
\renewcommand*{\Affilfont}{\small\it}
\date{} \maketitle
\vspace{-20pt}

\footnotetext{\hspace{-1.9mm}$^\dag$Corresponding author.\\
Email addresses: 20231577@dlnu.edu.cn(Q. Li); wolfzws@163.com(W. Zhou); yzhan@jlu.edu.cn(Y. Han).

\thanks{
$^*$The first author is supported by the National Natural Science Foundation of China (Grant No.~12501255) and  the Liaoning Provincial Natural Science Foundation of China (Grant No.~2025BS0268). The second author is supported by the National Natural Science Foundation of China (Grant No.~12571230).}}

{\bf Abstract}
This paper is  devoted to the study of normalized solutions to the following Kirchhoff type equation with  a logarithmic perturbation
\[
- \left( a + b \int_{\mathbb{R}^3} |\nabla u|^2 \,\mathrm{d}x \right) \Delta u  = \lambda u+ |u|^{p-2}u + u \log u^2, \quad x \in \mathbb{R}^3,
\]
under the normalized constraint $\int_{\mathbb{R}^3} u^2 \,\mathrm{d}x = c^2$, where $a,b>0$, $2<p\leq 6$, $c>0$ is a constant, and $\lambda\in\mathbb{R}$ emerges as a Lagrange multiplier  which is not a priori known.
A unified variational framework is developed based on Orlicz spaces together with  the Pohozaev constraint method and refined fiber map analysis.  For $2<p<\frac{14}{3}$ or $p=\frac{14}{3}$ with small mass, the energy functional is bounded from below and admits a positive radial ground state minimizer. For $\frac{14}{3}<p<6$, where the energy functional is unbounded from below, we establish the existence of two normalized solutions for 
 small mass: a ground state $u_c^+$ obtained via local minimization, and a second solution $u_c^-$ obtained via minimization on the negative component of the Pohozaev manifold. For the Sobolev critical case $p=6$, we construct a ground state solution  and, under a technical condition on the parameters, a second solution by introducing a proper auxiliary functional and precise energy estimates with Aubin-Talenti bubbles.
Asymptotically as $c\to0^+$,  the $L^{2}$ norm of the gradient of ground state solution  vanishes for  $2<p\le6$. Surprisingly, for $\frac{14}{3}<p<6$,  the $L^{2}$ norm of the gradient of the second solution diverges to infinity as $c\to 0^+$,  while for $p=6$ it concentrates around the  Aubin-Talenti bubble  with energy converging to the energy level of the corresponding critical Kirchhoff equation. Comparing the results in this paper with the
known ones, one sees that some new phenomena occur when the logarithmic perturbation is introduced.

Our work extends Soave's theory of normalized solutions  (\textit{J. Differential Equations} 2020 \& \textit{J. Funct. Anal.} 2020) to the critical Kirchhoff problems and provides a framework for a logarithmic perturbation in both the $L^2$-critical and   Sobolev critical cases.



\section{Introduction}
\setcounter{equation}{0}
This paper is devoted to the study of normalized solutions to the following Kirchhoff type equation with a logarithmic perturbation
\begin{equation}\label{eq:main-1}
- \left( a + b \int_{\mathbb{R}^3} |\nabla u|^2 \,\mathrm{d}x \right) \Delta u = \lambda u+ |u|^{p-2}u +  u \log u^2, \quad x \in \mathbb{R}^3,
\end{equation}
under the mass constraint
\begin{equation}\label{eq:main-2}
\int_{\mathbb{R}^3} u^2 \,\mathrm{d}x = c^2,
\end{equation}
where $c > 0$ is a constant, $a > 0$, $b > 0$, $p \in (2, 6]$ and $\lambda \in \mathbb{R}$ appears as a Lagrange multiplier.

For the one-dimensional case, equation \eqref{eq:main-1} is closely related to the
stationary analogue of the Kirchhoff type wave equation
\begin{eqnarray}\label{P2}
u_{tt}-(a+b\int_{\Omega}|\nabla u|^2\mathrm{d}x)\Delta u=f(x,u),
\end{eqnarray}
which was first proposed by Kirchhoff \cite{Kirchhoff} in 1883 to describe the transversal oscillations of a stretched string. In this mechanical model, the coefficients carry clear physical interpretations:
$u$ denotes the displacement of the string, $f$ stands for the external force,
$a$ represents the initial tension and $b$ is related to the intrinsic properties of the string.
One can refer to \cite{Arosio, Kirchhoff equations Yuzhu Han,He} and the references therein for more mathematical and physical background of \eqref{P2}.

Problem \eqref{eq:main-1} is usually referred to as being nonlocal,
since the presence of the term $(\int_{\mathbb{R}^3} |\nabla u|^2 \, \,\mathrm{d}x) \Delta u$
means that the equation is no longer a point-wise identity. When one considers the
existence of  solutions to problem \eqref{eq:main-1}, the main difficulty caused by
this nonlocal term is that generally one can not deduce from $u_n\rightharpoonup u$
weakly in $H^{1}(\mathbb{R}^3)$ the convergence $(\int_{\mathbb{R}^3} |\nabla u_{n}|^2 \,\mathrm{d}x) \Delta u_n\rightharpoonup (\int_{\mathbb{R}^3} |\nabla u|^2 \,\mathrm{d}x) \Delta u$ in $H^{-1}(\mathbb{R}^3)$.

In the classical variational problems with fixed $\lambda$, we point out that the existence, multiplicity and concentration of solutions for Kirchhoff type equation involving subcritical, critical and supercritical exponents have been extensively studied under different assumptions about the nonlinearity $f$. We refer the
interested reader to \cite{AlvesCorreaMa2005,Silva2019,ZhangDu2020}
and the references therein.

In the past few years, fruitful results have been obtained on the study of normalized solutions to elliptic equations. Such solutions are physically more meaningful because they can better reflect the dynamical properties of stationary solutions, including stability and instability \cite{Berestycki, Cazenave}.
In particular,  the work of Soave \cite{Soave2020a, Soave2020b} on Schr\"{o}dinger equations with combined nonlinearities  provides a framework for  studying the existence of normalized solutions.   In \cite{Ye2014,Ye2015,Ye2016}, Ye studied the problem
\begin{align}\label{eq:K1}
-\left(a + b \int_{\mathbb{R}^N} |\nabla u|^2  \,\mathrm{d}x \right) \Delta u = \lambda u + g(u) \quad \text{in } \mathbb{R}^N,
\end{align}
 under the normalized constraint $\int_{\mathbb{R}^N} u^2 \,\mathrm{d}x = c^2$ with $g(u)=|u|^{q-2}u$,
where 
$a > 0$, $b > 0$, $N\le3 $ and $2<q<2^{*}=:\frac{2N}{N-2}\ (2^{*}=\infty$  if $N=1,2)$. Ye's work systematically investigated the existence of weak solutions to \eqref{eq:K1} and the mass concentration behavior. 
For the $L^2$-subcritical case $2 < q < 2 + \frac{4}{N}$, a global minimizer exists for any  $c > 0$. For the intermediate range $2 + \frac{4}{N} \le q < 2 + \frac{8}{N}$, there exists a global minimizer  when  $c > c_q$. For the $L^2$-critical case $q = 2 + \frac{8}{N}$, no solution exists for  $0 < c \le c^*$, and a mountain pass type solution is obtained for $ c > c^*$. For the $L^2$-supercritical case where $2 + \frac{8}{N} < q < 2^*$, a mountain pass type solution  exists for all $c > 0$. Later, He et al. \cite{He2023} studied the existence and  blow-up behavior of  normalized solutions to  problem \eqref{eq:K1} with general nonlinearities in the $L^2$-subcritical case. Under suitable assumptions on $g(u)$, a ground state normalized solution  for any $c>0$ was obtained.

The study of normalized solutions to Kirchhoff equations, particularly those involving  the Sobolev critical exponent, has emerged as a significant research direction in nonlinear analysis. An interesting piece of work in this area was established by Zhang and Han \cite{Zhang2022}.
They  studied  the existence of
normalized solutions to the problem
\begin{equation}\label{eq:K}
-\Big(a+b\int_{\mathbb{R}^3}|\nabla u|^2\,\mathrm{d}x \Big)\Delta u = \lambda u + |u|^{p-2}u + \mu|u|^{q-2}u,
\end{equation}
under the normalized constraint \eqref{eq:main-2},
where $\mu = 1$,  $\frac{14}{3}\le q < 6$ and $p=6$,
and proved that the problem \eqref{eq:K}  has a  solution for all $c>0$.
Li et al. \cite{LiLuoYang2022} considered the existence and asymptotic properties of the solutions to the  Kirchhoff problem
\eqref{eq:K}--\eqref{eq:main-2} with $\mu>0$.
Their work  treated the mixed critical case $2<q<\frac{10}{3}$, $\frac{14}{3}<p\le 6$ and the purely $L^2$-supercritical case $\frac{14}{3}<q<p\le 6$, and established a multiplicity theory encompassing both local minimizers and mountain pass solutions.
Most recently, Chen and Tang \cite{ChenTang2025}  established a comprehensive existence theory for normalized solutions to  problem \eqref{eq:K} with $2<q<6$, $p=6$ and $\mu\in\mathbb{R}$. They showed the existence of two solutions (a local minimizer and a mountain pass solution) for $2<q<\frac{10}{3}$, addressed the gap for $\frac{10}{3}\le q<\frac{14}{3}$ by constructing a mountain pass solution, and obtained ground states for $\frac{14}{3}\le q<6$, along with the  nonexistence result for $\mu\le0$.
They introduced an auxiliary functional to overcome the difficulties arising from the nonlocal term $\|\nabla u\|_2^4$ and established the compactness thresholds through refined energy estimates.

A natural and important extension of this research line involves the inclusion of logarithmic perturbation terms, which arise in various physical contexts including quantum optics, nuclear physics, and Bose-Einstein condensation~\cite{AlfaroCarles2017}.
The main difficulty caused by the logarithmic nonlinearity  is the lack of smoothness of the associated energy functional.
 To overcome this difficulty, Cazenave~\cite{Cazenave1983} introduced the Banach space
\[
 \widetilde{W}:=\Big\{u\in H^1(\mathbb{R}^N)\,\Big|\,\int_{\mathbb{R}^N}u^2|\log u^2|\,\mathrm{d}x<\infty\Big\},
\]
equipped with an appropriate Luxemburg norm, which has become the standard framework for studying such problems. Following this approach, Shuai and Yang~\cite{ShuaiYang2025} investigated normalized solutions for the logarithmic Schrödinger equation with power  nonlinearity
\[
-\Delta u + \lambda u = \alpha u\log u^2 + \mu|u|^{p-2}u,
\]
under the mass constraint $\int_{\mathbb{R}^N}u^2\,\mathrm{d}x=c^2$,  where $ N \geq 2$, $\alpha, \mu\in \mathbb{R}$ and $2<p < 2^{*}=:\frac{2N}{N-2}\ (2^{*}=\infty$  if $N=1,2)$. Their work sheds some light on handling the interplay between logarithmic and power-type nonlinearities in the normalized setting.
More recently, Deng et al.
\cite{DengShiYang2025} studied the Choquard equation with logarithmic perturbation
\[
-\Delta u + \lambda u = u\log u^2 + \mu(I_\alpha\ast|u|^p)|u|^{p-2}u,
\]
under the mass constraint $\int_{\mathbb{R}^N}u^2\,\mathrm{d}x=c^2$, where $I_\alpha$ denotes the Riesz potential. Their work addresses the full range of critical exponents, including the Hardy-Littlewood-Sobolev upper critical exponent $p=\frac{N+\alpha}{N-2}$, the lower critical exponent $p=\frac{N+\alpha}{N}$, and the $L^2$-critical exponent $p=\frac{N+\alpha+2}{N}$. Through the  refined analysis of the Pohozaev manifold and careful energy estimates, they established the  existence of both ground state solutions and mountain pass solutions.

Inspired mainly by the above-mentioned results, our objective is to investigate normalized solutions to the logarithmic Kirchhoff problem \eqref{eq:main-1}--\eqref{eq:main-2}.  We develop a unified variational framework based on Orlicz spaces together with the Pohozaev constraint method and refined fiber map analysis.
For $2<p<\frac{14}{3}$ or $p=\frac{14}{3}$ with small mass, the functional energy is bounded from below, and it admits a global minimizer, which is a positive, radially symmetric ground state solution to \eqref{eq:main-1}--\eqref{eq:main-2}. For the mass-supercritical case $\frac{14}{3}<p<6$, the energy functional is unbounded from below. Using the Pohozaev constraint method, we showed the existence of two normalized solutions when $0<c< c_{**}$ for some $c_{**}>0$.
The Sobolev critical case $p=6$ brings some additional challenges due to the loss of compactness of the embedding $H^1(\mathbb{R}^3)\hookrightarrow L^6(\mathbb{R}^3)$. We prove the existence of a ground state solution $u_c^+$ via constrained minimization. Moreover, by introducing a proper auxiliary functional $\mathcal{J}$ to control the mountain pass level of the energy functional of problem \eqref{eq:main-1}--\eqref{eq:main-2} to satisfy the local compactness condition, we obtain a second solution $u_c^-$ (which is a mountain pass solution), under some technical conditions on the parameters $a, b$ and the optimal Sobolev constant $\mathcal{S}$. The construction of the second solution relies heavily on precise estimates on the Aubin-Talenti bubble which helps us to prove the local compactness  of the Palais-Smale sequences.
 Our asymptotic analysis reveals that as $c\to0^+$, the $L^{2}$ norm of the gradient of ground states  vanishes for  $p\le6$. Surprisingly, for $\frac{14}{3}<p<6$,  the $L^{2}$ norm of the gradient of the second solution diverges to infinity as $c\to 0^+$, while for $p=6$ it concentrates around the Aubin-Talenti bubble with its energy converging to the energy level of the corresponding critical Kirchhoff equation.

Comparing our results with that in \cite{Ye2014,Ye2015,Ye2016}, one sees that the logarithmic term plays a positive role for problem  to admit weak solutions  and this leads to some interesting phenomena.
This work extends Soave's theory \cite{Soave2020a, Soave2020b}  of normalized solutions from Schr\"{o}dinger equations to the critical Kirchhoff framework and provides the first
comprehensive research of logarithmic perturbations at both $L^2$-critical and Sobolev critical exponents.  

The organization of this paper is as follows.
We list some notations, definitions and necessary lemmas in Section 2. The main results of this paper are also stated here. 
In Section 3, we deal with the case $2<p\le \frac{14}{3}$ and prove Theorem \ref{thm:main1}. 
Some basic properties of the Pohozaev manifold and fiber maps are proved in Section 4 for the 
case $\frac{14}{3} < p \le 6$. In Section 5, we study the case $\frac{14}{3}<p< 6$ and prove 
Theorem  \ref{thm:main2}. In Section 6, we study the case $p=6$ and prove Theorems \ref{thm:main3} 
and \ref{thm:main4}. The asymptotic behavior of solutions are investigated in Section 7.

\section{Preliminaries and the main results}

In order to state our main results precisely, we first introduce some notations and
definitions and give some basic properties.
Throughout this paper, we adopt the following conventions. For $1 \le p \le \infty$, $\|\cdot\|_p$ denotes the standard $L^p(\mathbb{R}^3)$ norm. The Hilbert space $H^1(\mathbb{R}^3)$ is equipped with the norm $\|u\|_{H^1}^2 := \|u\|_2^2 + \|\nabla u\|_2^2$,  
and its dual space is written as $H^{-1}(\mathbb{R}^3)$. For each Banach space $X$, we use $\rightharpoonup$ and $\to$ to denote the weak and strong convergence in it, respectively. The letters $C, C_1, C_2, \ldots$ represent positive constants that may vary from line to line, while $B_R(x_0)$ denotes the open ball of radius $R$ centered at $x_0$ in $\mathbb{R}^3$. Let $\mathcal{S} > 0$ be the optimal Sobolev constant for the embedding $H^1(\mathbb{R}^3) \hookrightarrow L^6(\mathbb{R}^3)$, which satisfies $\mathcal{S}\|u\|_6^2 \le \|\nabla u\|_2^2$ for all $u \in H^1(\mathbb{R}^3)$. $O(t)$ denotes a quantity with $|\frac{O(t)}{t}| \le C$, $o(t)$ means $|\frac{o(t)}{t}| \to 0$ as $t \to 0$, and $o_n(1)$ denotes an infinitesimal as $n \to \infty$.

\subsection{Preliminaries}

One of the main  challenges in studying \eqref{eq:main-1} arises from the logarithmic perturbation  $u\log u^2$. 
Unlike the power-type nonlinearities, this term renders the energy functional not well-defined on $H^1(\mathbb{R}^3)$. 
Indeed, there exists $u \in H^1(\mathbb{R}^3)$ such that $\int_{\mathbb{R}^3} u^2 \log u^2 \,{\rm d}x = -\infty$. 
To overcome this difficulty, we employ the Orlicz space framework originally developed by Cazenave \cite{Cazenave1983} 
for logarithmic Schr\"{o}dinger equations.

We begin by introducing the convex functions $A, B: [0,\infty) \to \mathbb{R}$:
\begin{equation*}
A(s) = \begin{cases}
-s^2 \log s^2, & 0 \le s < e^{-3}, \\
3s^2 + 4e^{-3}s - e^{-6}, & s \ge e^{-3},
\end{cases}
\qquad
B(s) = s^2 \log s^2 + A(s).
\end{equation*}

The following lemma summarizes some basic properties of $A(s)$ and $B(s)$ (see \cite{Cazenave1983,LiebLoss2001,Shuai2019} for the proof).

\begin{lemma}\label{A-B}
The functions $A$ and $B$ defined above satisfy:
\begin{itemize}
\item[(i)] $A$ and $B$ are positive, convex, and increasing on $(0,\infty)$.
\item[(ii)] There exists $K_q > 0$ such that  $B(s) \le K_q s^q$ for all $q \in (2, \frac{10}{3}) $ and $s > 0$.
\item[(iii)]  There exists $C>0$ such that  $A(x+ty)-A(x) \le CtA(y)$ for all $x,y,t\geq 0$ with
$ ty\ge \frac{x}{2}$.
\end{itemize}
\end{lemma}

The Orlicz space $V$ associated with $A$ is defined as
\begin{equation*}
V := \left\{ u \in L^1_{\text{loc}}(\mathbb{R}^3) \,\bigg|\, A(|u|) \in L^1(\mathbb{R}^3) \right\},
\end{equation*}
endowed with the Luxemburg norm
\begin{equation*}\label{eq:Luxemburg}
\|u\|_V := \inf\left\{ k > 0 \,\bigg|\, \int_{\mathbb{R}^3} A\left(\frac{|u|}{k}\right) \,\mathrm{d}x \le 1 \right\}.
\end{equation*}

As for the structural properties of $V$, we have the following lemma.

\begin{lemma}\label{lem:Orlicz-V}
The space $V$ defined above satisfies
\begin{itemize}
\item[(i)] $V$ is a reflexive Banach space.
\item[(ii)] For any $u \in V$, we have
\[
\min\{\|u\|_V, \|u\|_V^2\} \le \int_{\mathbb{R}^3} A(|u|)
\,\mathrm{d}x \le \max\{\|u\|_V, \|u\|_V^2\}.
\]
\item[(iii)] If $u_n \to u$ a.e. in $\mathbb{R}^3$ and $\int_{\mathbb{R}^3} A(|u_n|) \,\mathrm{d}x \to \int_{\mathbb{R}^3} A(|u|) \,\mathrm{d}x < \infty$, then $\|u_n - u\|_V \to 0$.
\end{itemize}
\end{lemma}

\begin{proof}
These are standard results in Orlicz space theory, and we refer to \cite{Cazenave1983} for a detailed exposition.
\end{proof}

We now define our working variational space as
\begin{equation*}\label{eq:W-space}
W := \left\{ u \in H^1(\mathbb{R}^3) \,\bigg|\,  \int_{\mathbb{R}^3}u^2|\log u^2|\,\mathrm{d}x<\infty\right\},
\end{equation*}
equipped with the norm
\[
\|u\|_W := \|u\|_{H^1} + \|u\|_V.
\]

The radial subspace $W_r := W \cap H^1_r(\mathbb{R}^3)$, where $H^1_r(\mathbb{R}^3)$  
denotes radial functions in $H^1(\mathbb{R}^3)$, will play a crucial role in our variational arguments. The following compactness result is fundamental; see \cite[Proposition 3.1]{Cazenave1983}.

\begin{lemma}\label{lem:compact-embedding}
The embedding $W_r \hookrightarrow L^2(\mathbb{R}^3)$ is compact.
\end{lemma}

\begin{lemma}[Gagliardo-Nirenberg  inequality \cite{Weinstein1983}]\label{lem:G-N}
Let $N \geq 1$ and $p \in (2, 2^*)$, where $2^* := \frac{2N}{N-2}$ for $N \geq 3$ is the Sobolev critical exponent
(and $2^* = \infty$ when $N=1,2$). Then for any $u \in H^1(\mathbb{R}^N)$, the following Gagliardo-Nirenberg inequality holds
\begin{equation*}
\|u\|_{L^p(\mathbb{R}^N)} \leq C_{N,p} \|\nabla u\|_{L^2(\mathbb{R}^N)}^{\gamma_p} \|u\|_{L^2(\mathbb{R}^N)}^{1-\gamma_p},
\end{equation*}
where $C_{N,p} > 0$ is a constant depending only on $N$ and $p$, and $\gamma_p = \frac{N(p-2)}{2p}$ is the interpolation exponent.
\end{lemma}

\begin{lemma}(\cite{Shuai2019})\label{lem;duishu}
Let $\{u_n\}$ be a bounded sequence in $H^1(\mathbb{R}^N)$ such that $u_n \to u$ a.e. in $\mathbb{R}^N$ and $\{u_n^2 \log u_n^2\}$ is a bounded sequence in $L^1(\mathbb{R}^N)$. Then $u^2 \log u^2 \in L^1(\mathbb{R}^N)$ and
\begin{equation*}\label{eq1-
lem;duishu}
\lim_{n\to\infty} \int_{\mathbb{R}^N} \bigl(u_n^2 \log u_n^2 - |u_n - u|^2 \log |u_n - u|^2\bigr)\,\mathrm{d}x = \int_{\mathbb{R}^N} u^2 \log u^2 \,\mathrm{d}x.
\end{equation*}
\end{lemma}

Solutions to problem \eqref{eq:main-1}--\eqref{eq:main-2} can be obtained by finding critical points of the energy functional
\begin{equation*}\label{eq:energy-functional}
I(u) := \frac{a}{2}\int_{\mathbb{R}^3} |\nabla u|^2 \,\mathrm{d}x
 + \frac{b}{4}\left(\int_{\mathbb{R}^3} |\nabla u|^2 \,\mathrm{d}x
\right)^2 - \frac{1}{p}\int_{\mathbb{R}^3} |u|^p \,\mathrm{d}x
 + \frac{1}{2}c^2 - \frac{1}{2}\int_{\mathbb{R}^3} u^2 \log u^2 \,\mathrm{d}x
,
\end{equation*}
under the constraint
\[
S(c):=\left\{u\in W \mid \|u\|_2=c\right\}.
\]

The $L^2$-preserving scaling $s * u(x) = e^{\frac{3s}{2}} u(e^s x)$, $s\in \mathbb{R}$ plays a central role in our analysis.
Then the fiber map $\Phi_u(s) := I(s *  u)$ has the following form
\begin{equation*}\label{eq:fiber-map}
\Phi_u(s) := \frac{a e^{2s}}{2}\|\nabla u\|_2^2 + \frac{b e^{4s}}{4}\|\nabla u\|_2^4 - \frac{e^{p\gamma_p s}}{p}\| u\|_p^p + \frac{c^2}{2}(1-3s) - \frac{1}{2}\int_{\mathbb{R}^3} u^2 \log u^2\,\mathrm{d}x,
\end{equation*}
where $\gamma_p = \frac{3(p-2)}{2p}$.

We define the Pohozaev functional
\begin{equation*}\label{eq:pohozaev}
Q(u) :=\Phi'_u(0)= a\|\nabla u\|_2^2  + b\|\nabla u\|_2^4 - \gamma_p\| u\|_p^p - \frac{3 c^2}{2},
\end{equation*}
which characterizes the critical points of $\Phi_u$. The associated Pohozaev manifold is
$$\Lambda_c := \{u \in S(c) \mid Q(u) = 0\}.$$
It is well known that any critical point of
$I|_{S(c)}$ stays in $\Lambda_c$, as a consequence of the Pohozaev identity (see \cite[Lemma 2.7]{Jeanjean1997}).
We further decompose $\Lambda_c$ into three disjoint subsets based on the second variation of the fiber map
\begin{align*}
\Lambda_c^0 &:= \{u \in \Lambda_c \mid \Phi_u''(0) = 0\}, \\
\Lambda_c^+&:= \{u \in \Lambda_c \mid \Phi_u''(0) > 0\}, \\
\Lambda_c^-&:= \{u \in \Lambda_c \mid \Phi_u''(0) < 0\}.
\end{align*}

We present several technical results that underpin our variational approach.

\begin{lemma}\label{lem:tangent}(\cite{BartschSoave2019})
For $u \in S(c)$, the tangent space to the constraint manifold $S(c)$ at $u$ is
\[
T_u S(c) = \left\{ v \in H^1(\mathbb{R}^3) \,\bigg|\, \int_{\mathbb{R}^3} uv \,\mathrm{d}x = 0 \right\}.
\]
Moreover, for any $s \in \mathbb{R}$, the scaling map $v \mapsto s*v$ induces a linear isomorphism $T_u S(c) \to T_{s * u} S(c)$ with inverse $\phi \mapsto (-s)*\phi$.
\end{lemma}

\begin{lemma}\label{lem:radial-decay}(\cite{Berestycki1983})
Let $u \in L^t(\mathbb{R}^3), t\in[1,\infty)$ be nonnegative,  radial and non-increasing (i.e., $0\le u(y)\le u(x)$ if $|x|\le |y|$) for $1 \le t < \infty$. Then
\[
|u(x)| \le |x|^{-\frac{3}{t} }\left( \frac{|\mathbb{S}^2|}{3} \right)^{-\frac{1}{t}} \|u\|_t, \quad x \neq 0,
\]
where $|\mathbb{S}^2|$ denotes the surface area of the unit sphere in $\mathbb{R}^3$.
\end{lemma}

\begin{definition}\label{def1}
We say that $u_0$ is a ground state solution to problem  \eqref{eq:main-1}--\eqref{eq:main-2} on $S(c)$ if
\[
\,\mathrm{d}I|_{S(c)}(u_0)=0\quad\text{and}\quad I(u_0)=\inf\left\{I(u)\mid \,\mathrm{d}I|_{S(c)}(u)=0\text{ and }u\in S(c)\right\}.
\]
\end{definition}
\begin{definition}\label{def:homotopy-stable}(\cite[]{Ghoussoub})
Let $X$ be a topological space and $\partial\Gamma$ be a closed subset of $X$.
We  say that a class $\Gamma$ of compact subsets of $X$ is a  homotopy stable family with extended boundary $\partial\Gamma$ if for any set $G$ in $\Gamma$ and any $\eta \in C([0, 1] \times X; X)$ satisfying $\eta(t, x) = x$ for all $(t, x) \in (\{0\} \times X) \cup ([0, 1] \times \partial\Gamma)$ we have that $\eta(\{1\} \times G) \in \Gamma$.
\end{definition}

We recall the min-max principle of Ghoussoub~\cite[Theorem~5.2]{Ghoussoub} in the following lemma.

\begin{lemma}\label{min-max principle}
Let $\Psi$ be a $C^1$-functional on a complete connected $C^1$-Finsler manifold $X$ and consider a homotopy-stable family
$\Gamma$ with an extended closed boundary
$\partial \Gamma$. Set
$$m = m(\Psi, \Gamma):=\inf\limits_{G \in \Gamma} \max\limits_{x \in G} \Psi(x)$$ and  $\sup \Psi(\partial\Gamma) < m $.
Then for any sequence of sets $\{\widetilde{G}_{n}\}$ in $\Gamma$ such that $\lim\limits_{n\to+\infty} \sup\limits_{\widetilde{G}_{n}} \Psi = m$, there exists a sequence $\{x_n\}$ in $X$ such that
\[
\lim_{n\to+\infty} \Psi(x_n) = m,\quad \lim_{n\to+\infty}  \|\mathrm{d}\Psi(x_n)\|= 0,\quad 
\lim_{n\to+\infty} \operatorname{dist}(x_n, \widetilde{G}_{n}) = 0.
\]
\end{lemma}

\begin{lemma}\label{lem:st}
There holds
\begin{enumerate}
\item[(i)] $(1+s)^{\frac{3}{2}}\le 1+2s$, for  sufficiently small $s>0$.
\item[(ii)] $(1+st)^{\frac{3}{2}}-1\le t^{\frac{3}{2}}[(1+s)^{\frac{3}{2}}-1]$ for any $s\ge 0$, $t\ge 1$.
\item[(iii)] $(1+s+t)^{\frac{3}{2}}-1\ge [(1+s)^{\frac{3}{2}}-1]+[(1+t)^{\frac{3}{2}}-1]$ for any $s$, $t\ge 0$.
 \item[(iv)] For any $\alpha > 0$, there exists $C_{\alpha}>0$ such that $\log t\le C_{\alpha}t^{\alpha}$ for $t>0$.
\end{enumerate}

\end{lemma}

\begin{lemma}\label{lem:fmax}
Assume that $Y>0$ is a constant. Set $f(t):=\mathcal{S}^{\frac{3}{2}} \Big[ \frac{a+bY}{2} t^2 + \frac{b \mathcal{S}^{\frac{3}{2}}}{4} t^4 - \frac{1}{6} t^6\Big] $ for $t>0$. Then
 \begin{align}\label{Gai-1}
f(t)\le C(a,b, \mathcal{S}) + \frac{\left(b^2 \mathcal{S}^4 + 4a\mathcal{S}\right)^{\frac{3}{2}}}{24} \left[ \left(1 + \frac{4b}{b^2 \mathcal{S}^3 + 4a} Y\right)^{\frac{3}{2}} - 1 \right] + \frac{b^2 S^3}{4} Y,
\end{align}
where $$C(a,b, \mathcal{S}):=\frac{ab\mathcal{S}^3}{4} + \frac{b^3 \mathcal{S}^6}{24} + \frac{(b^2 \mathcal{S}^4 + 4a\mathcal{S})^{3/2}}{24}.$$
\end{lemma}

\begin{proof}
It is directly verified that $f(0)=0$, $f(t)>0$ for small $t$, $f(t) \to -\infty$ as $t \to +\infty$, and $f(t)$ attains its maximum at $t_0$ with
 $$t_0^2= \frac{1}{2}\left[b \mathcal{S}^{\frac{3}{2}} + \sqrt{b^2 \mathcal{S}^3 + 4\left(a + bY\right)}\right].$$
The right hand side of \eqref{Gai-1} is exactly $f(t_0)$.
The proof is complete.
\end{proof}

\subsection{The main results}

In this section, we state our results concerning normalized solutions to problem \eqref{eq:main-1}--\eqref{eq:main-2}. The results are categorized according to the range of the exponent $p\in(2,6]$, revealing the distinct geometric structures of the fiber map.

We first consider the case $2<p\leq \frac{14}{3}$ where the energy functional is bounded from below on $S(c)$, 
and establish the existence of the global minimizer. For this, set
\begin{align}\label{gai-5}
m(c) := \inf_{u\in S(c)} I(u).
\end{align}

\begin{theorem}\label{thm:main1}
Let $p$  satisfy either
\begin{itemize}
    \item[(i)] $2<p<\frac{14}{3}$, or
    \item[(ii)] $p=\frac{14}{3}$ with $0<c< \widetilde{c}$, where $\widetilde{c}$ is defined in Lemma \ref{lem:m(c)-bounded-below}.
\end{itemize}
Then $m(c)$ is attained at a bounded, positive and  radially symmetric  function $u_c\in S(c)$.
Moreover, there exists  a $\lambda_c\in\mathbb{R}$ such that $(\lambda_c,u_c)$ solves  problem \eqref{eq:main-1}--\eqref{eq:main-2} and $u_c$ is a ground state solution. \end{theorem}

For $\frac{14}{3}<p\le6$, the energy functional is no longer bounded from  below on $S(c)$. However, by employing the Pohozaev constraint method, we uncover a much richer geometric structure of the fiber map, yielding two normalized solutions for small mass.
Set
\begin{align}\label{gai-6}
m_c^\pm := \inf_{u\in\Lambda_c^\pm} I(u).
\end{align}

\begin{theorem}
\label{thm:main2}
Let $\frac{14}{3}<p<6$ and $c\in(0,c_{**})$, where $c_{**}$ is the constant determined in Lemma \ref{lem:critical_mass}.
Then there exist two couples of solutions $(\lambda_c^\pm,u_c^\pm)\in\mathbb{R}\times S(c)$ solving problem \eqref{eq:main-1}--\eqref{eq:main-2} with the following properties.
\begin{itemize}
    \item[(i)] $u_c^+\in\Lambda_c^+$ is a local minimizer of $I(u)$ on $S(c)$ with $I(u_c^+)=m_c^+$.
    \item[(ii)] $u_c^-\in\Lambda_c^-$ is a mountain pass  solution with $I(u_c^-)=m_c^-$.
    \item[(iii)] Both $u_c^\pm$ are bounded, positive,  and radially symmetric.
    \item[(iv)] $u_c^+$ is a ground state solution.\end{itemize}
\end{theorem}

The Sobolev critical case $p=6$ presents some additional challenges due to the loss of compactness
of the embedding $H^1(\mathbb{R}^3)\hookrightarrow L^6(\mathbb{R}^3)$.
 The following result establishes the existence of a ground state solution.

\begin{theorem}
\label{thm:main3}Let $p=6$ and $c\in(0,\widetilde{c}_{**})$, where $\widetilde{c}_{**}$ is defined in Lemma~\ref{lem:critical_mass}. Then
there exists a bounded, positive, and radially symmetric  solution $u_c^+\in\Lambda_c^+$ satisfying $I(u_c^+)=m_c^{+}$.
Moreover, there exists a $\lambda_c^+\in\mathbb{R}$ such that $(\lambda_c^+,u_c^+)$ solves  problem \eqref{eq:main-1}--\eqref{eq:main-2} and $u_c^+$ is a ground state solution.

\end{theorem}

To obtain a second solution, we introduce an auxiliary functional $\mathcal{J}(u)$ defined by \eqref{gai-J}, 
which plays a crucial role  in controlling the mountain pass level of the energy functional associated with problem \eqref{eq:main-1}--\eqref{eq:main-2} from above to satisfy the local compactness condition.

 \begin{theorem}
\label{thm:main4}
Let $p=6$ and $c\in(0,\min\{\widetilde{c}_{**},c_0,c_1,c_{2}\})$, where $\widetilde{c}_{**}$, $c_0$, $c_2$ and $c_1$ are defined  in Lemmas~\ref{lem:critical_mass},
\ref{lem:c_02},  \ref{lem2:c0c1} and \eqref{c11}, respectively. Assume that the technical condition \eqref{eq:technical} holds. Then there exists a bounded, positive, and radially symmetric mountain pass solution $u_c^-\in\Lambda_c^-$ satisfying $I(u_c^-)=m_c^-$. Moreover, there exists a $\lambda_c^-\in\mathbb{R}$ such that $(\lambda_c^-,u_c^-)$ solves  problem \eqref{eq:main-1}--\eqref{eq:main-2}.
\end{theorem}

\begin{remark}
 Comparing our results with the ones obtained by Ye in \cite{Ye2014,Ye2015,Ye2016}, one can see that the logarithmic perturbation $u\log u^2$ substantially enhances the existence theory of such problem.

(i) For $ \frac{10}{3} < p \le \frac{14}{3}$, Ye established the existence result for certain $c>0$,
while we obtain one solution to  problem \eqref{eq:main-1}--\eqref{eq:main-2}, which is a global minimizer of the corresponding functional, for any $c > 0$.

(ii) For $\frac{14}{3} < p < 6$, only a mountain pass solution was known by Ye. Here, we establish the existence of two solutions for small mass: a ground state solution and a second solution, which induces the multiplicity of solutions.

(iii) Additionally, we consider the Sobolev critical case $p = 6$, and provides the existence of two solutions, which   completes the existence theory for the full range of $p\in(2,6]$.
\end{remark}

We conclude with a characterization of the asymptotic properties of solutions as $c\to 0^+$.

\begin{theorem}
\label{th:to}
Let $u_c^\pm$ be the solutions obtained in Theorems \ref{thm:main1}-\ref{thm:main4}. 
Then
\begin{itemize}
    \item[(i)] For all $2<p\leq 6$, the ground state solution satisfies $\|\nabla u_c^+\|_2^2\to 0$, $m_c^+\to 0$, and $m(c)\to 0$ as $c\to 0^+$.

    \item[(ii)] For $\frac{14}{3}<p<6$, the second solution satisfies $\|\nabla u_c^-\|_2^2 \to +\infty$ as $c\to 0^+$.

    \item[(iii)] For $p=6$, the second solution satisfies $\|\nabla u_c^-\|_2^2\to \frac{\mathcal{S}}{2}\Big(b\mathcal{S}^2 + \sqrt{  b^2\mathcal{S}^4+4 a\mathcal{S}}\Big)$ and $I(u_c^-)\to C(a,b,\mathcal{S})$ as $c\to 0^+$, where $ C(a,b,\mathcal{S})$  is defined in Lemma \ref{lem:fmax}.
\end{itemize}
\end{theorem}

\begin{remark}
The asymptotic behavior in (iii) reveals a striking phenomenon: while the first solution vanishes, the second solution concentrates around the Aubin-Talenti bubble with its energy converging to the energy level of the corresponding critical Kirchhoff equation.
\end{remark}

\section{ $L^2$-subcritical and $L^2$-critical cases: $2<p \leq \frac{14}{3}$}
\setcounter{equation}{0}

In this section, we study the existence of global minimizers for $I(u)$ on the constraint manifold $S(c)$ in both $L^2$-subcritical and $L^2$-critical cases. We first show that the energy functional is bounded from below on $S(c)$ under the appropriate conditions.

\begin{lemma}\label{lem:m(c)-bounded-below}
Suppose that either
\begin{itemize}
    \item[(i)] $2<p<\frac{14}{3}$, or
    \item[(ii)] $p=\frac{14}{3}$ with $0<c< \widetilde{c}=(\frac{bp}{4C^{p}_{3,p}})^{\frac{1}{(1-\gamma_p)p}}$.
\end{itemize}
Then $I(u)\rightarrow+\infty$ as $\|\nabla u\|_2\rightarrow+\infty$ and consequently $m(c) > -\infty$, where $m(c)$ is defined in \eqref{gai-5}.\end{lemma}

\begin{proof}
For any $u \in S(c)$, by the Gagliardo-Nirenberg inequality and property (ii) of $B(u)$,  we have
\begin{align}I(u) &= \frac{a}{2}\|\nabla u\|_2^2 + \frac{b}{4}\|\nabla u\|_2^4 - \frac{1}{p}\|u\|_p^p + \frac{1}{2}c^2 + \frac{1}{2}\int_{\mathbb{R}^3}\bigl(A(|u|) - B(|u|)\bigr)\,\mathrm{d}x\nonumber\\
&\geq \frac{a}{2}\|\nabla u\|_2^2 + \frac{b}{4}\|\nabla u\|_2^4 - \frac{1}{p}C^{p}_{3,p}\|\nabla u\|_2^{p\gamma_p}c^{p(1-\gamma_p)}
- \frac{ K_q C^{p}_{3,q}}{2}\|\nabla u\|_2^{q\gamma_q}c^{q(1-\gamma_q)}.\label{lem:m(c)-bounded-below-eq1}
\end{align}

In case (i), we have $p\gamma_p < 4$  and $q\gamma_q < 2$, which implies that the nonlocal term $\frac{b}{4}\|\nabla u\|_2^4$ dominates as $\|\nabla u\|_2 \to \infty$. 
 Consequently, $I(u)$ is coercive and bounded from below. 
In case (ii), we have $p\gamma_p = 4$ and the condition $c<\widetilde{c}$ ensures that the coefficient 
$\frac{b}{4} - \frac{1}{p}C^{p}_{3,p}c^{(1-\gamma_p)p}$ of the term $\|\nabla u\|_2^4$ is positive,
which again yields coercivity. This completes the proof.
\end{proof}

To prove the existence of solutions to problem \eqref{eq:main-1}--\eqref{eq:main-2}, we need an equivalent characterization of $m(c)$ via the radially symmetric functions by defining
\[ m_r(c) := \inf_{u \in S(c) \cap H_r^1(\mathbb{R}^3)} I(u).\]
 \begin{lemma}\label{lem:m(c)=mr(c)}
Under the assumptions of Lemma \ref{lem:m(c)-bounded-below}, we have $m(c) = m_r(c)$.
\end{lemma}

\begin{proof}
Since $S(c) \cap H_r^1(\mathbb{R}^3) \subset S(c)$, we immediately have $m(c) \leq m_r(c)$. For the reverse inequality, given any $u \in S(c)$, let $u^*$ denote its symmetric decreasing rearrangement, which  is a positive, radially symmetric
non-increasing function. By the Riesz rearrangement inequality \cite{LiebLoss2001}
 and the convexity of $A(u)$ and $B(u)$, we have
\[
\|\nabla u^*\|_2 \leq \|\nabla u\|_2, \quad \|u^*\|_l = \|u\|_l \quad \text{for} \quad l\in[2,6],\] and
\[\int_{\mathbb{R}^3}A(|u^*|)\,\mathrm{d}x = \int_{\mathbb{R}^3}A(|u|)\,\mathrm{d}x,\quad \int_{\mathbb{R}^3}B(|u^*|)\,\mathrm{d}x = \int_{\mathbb{R}^3}B(|u|)\,\mathrm{d}x.\]
Consequently, $I(u^*) \leq I(u)$. Since $u^* \in S(c) \cap H_r^1(\mathbb{R}^3)$, taking the infimum over all $u \in S(c)$ yields $m_r(c) \leq m(c)$. Therefore, $m(c) = m_r(c)$. The proof is complete.
\end{proof}
\begin{proof}[Proof of Theorem \ref{thm:main1}] From Lemma \ref{lem:m(c)=mr(c)}, we only need to prove that  there exists a function $u\in S(c) \cap H_r^1(\mathbb{R}^3)$ such that $I(u)=m_r(c)$.
Let $\{u_n\} \subset S(c) \cap H_r^1(\mathbb{R}^3)$ be a minimizing sequence for $m_r(c)$.
By Lemma \ref{lem:m(c)-bounded-below}, $\{u_n\}$ is bounded in $H^1(\mathbb{R}^3)$, 
and therefore $\{B(|u_n|)\}$ is bounded in $L^1(\mathbb{R}^3)$. Then it follows from \eqref{lem:m(c)-bounded-below-eq1}  that $\{A(|u_n|)\}$ is also bounded in $L^1(\mathbb{R}^3)$. The definition of the space $W$ and Lemma \ref{lem:Orlicz-V} imply that $\{u_n\}$ is also bounded in $W_{r}$. Therefore, up to a subsequence if necessary, we have
$$u_n \rightharpoonup u \quad \text{weakly in } W_{r}.$$ By Lemma~\ref{lem:compact-embedding}, we know that $u_n \to u$ strongly in $L^2(\mathbb{R}^3)$  and almost everywhere in $\mathbb{R}^3$. Consequently,
we know that $u_n \to u$ strongly in $L^l(\mathbb{R}^3)$ for all $l \in \left(2, 6\right)$ by Lemma \ref{lem:G-N}, 
which in turn implies that
$B(|u_n|) \to B(|u|) $ strongly in $L^1(\mathbb{R}^3)$.
By the weak lower semi-continuity of the functional and the fact that $u\in S(c)\cap H_r^1(\mathbb{R}^3)$, we further derive
\begin{align*}
   m_r(c)  &=\liminf_{n \to \infty} I(u_n)\\
   &=
\liminf_{n \to \infty}\left[\frac{a}{2} \|\nabla u_n\|_{2}^{2} + \frac{b}{4} \|\nabla u_n\|_{2}^{4}-  \frac{1}{p}  \|u_n\|_{p}^{p} + \frac{1}{2} c^{2}+\frac{1}{2} \int_{\mathbb{R}^3} \bigl(A(|u_n|) - B(|u_n|) \bigr) \,\mathrm{d}x\right] \\
& \geq
\liminf_{n \to \infty}\left(\frac{a}{2} \|\nabla u_n\|_{2}^{2} + \frac{b}{4} \|\nabla u_n\|_{2}^{4}+\frac{1}{2} \int_{\mathbb{R}^3} A(|u_n|) \,\mathrm{d}x\right) + \frac{1}{2} c^{2}- \frac{1}{p}  \|u\|_{p}^{p}- \frac{1}{2} \int_{\mathbb{R}^3} B(|u|) \,\mathrm{d}x  \\
  &\geq \frac{a}{2} \|\nabla u\|_{2}^{2} + \frac{b}{4} \|\nabla u\|_{2}^{4}-  \frac{1}{p}  \|u\|_{p}^{p} + \frac{1}{2} c^{2}+\frac{1}{2} \int_{\mathbb{R}^3} \left(A(|u|) - B(|u|)\right) \,\mathrm{d}x\\
  &
   = I(u)  \geq m_r(c).
\end{align*}
Therefore, we have
$I(u) = m_r(c)$, $u_n \to u$ strongly in $ H^1(\mathbb{R}^3)$ and $A(|u_n|) \to A(|u|)$ strongly in $ L^1(\mathbb{R}^3)$.
An application of Lemma~\ref{lem:Orlicz-V} ensures that \(u_n \to u\) strongly in \(V\) as \(n \to \infty\) and hence \(u_n \to u\) strongly in \(W_r\).

Since $I$ is even, we may assume that $u \geq 0$.
By the Lagrange multiplier rule (see \cite[Lemma 3]{Berestycki1983}), there exists a  $\lambda \in \mathbb{R}$ such that $u$  satisfies \eqref{eq:main-1}.
By the regularity theory of elliptic equations, 
we deduce that $u \in C^2(\mathbb{R}^3)$. Set $\Theta := a + b\int_{\mathbb{R}^3} |\nabla u|^2\,dx > 0$. 
Then equation \eqref{eq:main-1} can be rewritten as
\[
\Theta \Delta u =- \lambda u -  u\log |u|^2 - |u|^{p-2}u \quad \text{in } \mathbb{R}^3.
\]
 Define $\beta(s) = -\lambda s -  s\log s^2$ for $s > 0$, and extend it continuously by $\beta(0) = 0$. Then $\beta \in C([0,\infty))$ is non-decreasing for $s $ small, with $\beta(0) = 0$ and $\beta({\sqrt{e^\lambda}}) = 0$. Moreover, for $u > 0$ sufficiently small, we have
\[
\Theta \Delta u \leq \beta(u).
\]
Applying the strong maximum principle from \cite[Theorem 1]{Vazquez1984}, we conclude that $u > 0$ in $\mathbb{R}^3$. Therefore, $u \in S(c)$ is a positive, radially symmetric function that attains $m(c)$. Then the conclusion of Theorem \ref{thm:main1} follows by taking $(\lambda_c, u_c)=(\lambda, u)$. The proof of Theorem \ref{thm:main1} is complete.
\end{proof}

\section{Basic properties of the Pohozaev manifold and the fiber maps for $\frac{14}{3} < p \le 6$}

\setcounter{equation}{0}

When$\frac{14}{3} < p \le  6$, the energy functional $I$ is no longer bounded from below on the constraint manifold $S(c)$, making the direct minimization approach ineffective. To overcome this difficulty, this section conducts a detailed analysis of  the fiber map  and the Pohozaev manifold. This geometric framework will be crucial for constructing both the local minimizer and the mountain pass solution in the subsequent sections.

For any $c > 0$, we consider the function
\begin{equation}\label{Gai-2}
\tilde{g}_c(r):=\begin{cases}
ar^2 + br^4 - \gamma_p C_{3,p}^p c^{p(1-\gamma_p)}r^{p\gamma_p},  & \frac{14}{3} <p<6,\\[4pt]
ar^2 + br^4 - \mathcal{S} ^{-3} r^{6}, & p=6,
\end{cases}\ \ \ \text{for}\  r>0.
\end{equation}

\begin{lemma}\label{lem:critical_mass}

Assume that $\frac{14}{3} < p \leq 6$. Then 
 there exists a $c_* > 0$  such that for any $c \in (0,c_*)$, 
the equation
\[
\tilde{g}_c(r) = \frac{3c^2}{2}
\]
has exactly two positive solutions, denoted by $R_0$ and $R_1$. Let $R_*$ be the unique maximizer of $\tilde{g}_c$ on $(0,\infty)$. 
Then we have
$
0 < R_0 < R_* < R_1,
$
and $\tilde{g}_c(r) > \frac{3c^2}{2}$ holds if $r \in (R_0,R_1)$. 
\end{lemma}

\begin{proof}
We divide the proof into two cases based on the range of  of $p$.

\textbf{Case 1: $\frac{14}{3} <p<6$.}
For each fixed $c>0$, a direct calculation shows that
 $\tilde{g}_c(0) = 0$, $\tilde{g}_c(r)>0$ for small $r$, $\lim\limits_{r\to\infty}\tilde{g}_c(r) = -\infty$, and
 $\tilde{g}_c$ has a unique global maximizer $R_c > 0$, which is characterized by
\begin{equation}\label{eq:c_star1}
2aR_c+4bR_c^3=\gamma_p^2pC_{3,p}^p\,c^{p(1-\gamma_p)}R_c^{p\gamma_p-1}.
\end{equation}
 Then the corresponding maximum value is
\begin{align*}
\tilde{g}_c(R_c) = a R_c^2 + b R_c^4 - \gamma_p C_{3,p}^p c^{p(1-\gamma_p)}R_c^{p\gamma_p}.
\end{align*}

Consider the equation $\tilde{g}_c(R_c) =\frac{3 c^2}{2}$, i.e.,
\begin{equation}\label{eq:c_star2}
aR_c^2 + bR_c^4 - \gamma_p C_{3,p}^p c^{p(1-\gamma_p)} R_c^{p\gamma_p} = \frac{3c^2}{2}.
\end{equation}

Set $z:= R_c^{2} > 0$ and $\beta := p(1 - \gamma_p) >0$. Then equations \eqref{eq:c_star1} and \eqref{eq:c_star2} become  \begin{align}
&2a + 4b z =  \gamma_p^2 p C_{3,p}^p c^{\beta} z^{\frac{p\gamma_p}{2} - 1}.
\label{1a}\\
&a z + b z^2 -  \gamma_p C_{3,p}^p c^{\beta} z^{\frac{p\gamma_p}{2}} = \frac{3 c^2}{2},
\label{2a}
\end{align}
where $p\gamma_p > 4$.
It follows from \eqref{1a} that
\begin{align}
 \gamma_p C_{3,p}^p c^{\beta} z^{\frac{p\gamma_p}{2}} = \frac{2a}{p \gamma_p} z + \frac{4b}{p \gamma_p} z^2.
\label{3}
\end{align}
Substituting \eqref{3} into \eqref{2a}, we obtain
\[
a z + b z^2 - \left( \frac{2a}{p \gamma_p} z + \frac{4b}{p \gamma_p} z^2 \right) = \frac{3 c^2}{2}.
\]
This yields the equivalent  equations
\begin{align}
&2a z^{1 - \frac{p\gamma_p}{2} } + 4b z^{2 - \frac{p\gamma_p}{2} } =  \gamma_p^2 p C_{3,p}^p c^{\beta}, \label{eq:3} \\
&a \left(1 - \frac{2}{p\gamma_p}\right) z + b \left(1 - \frac{4}{p\gamma_p}\right) z^2 = \frac{3 c^2}{2}. \label{4}\nonumber
\end{align}
Now consider the function $H(z) := 2a z^{1 - \frac{p\gamma_p}{2}} + 4b z^{2 - \frac{p\gamma_p}{2}}$ for $z > 0$. 
Since $p\gamma_p > 4$, we have $1 - \frac{p\gamma_p}{2}< -1$ and $2 - \frac{p\gamma_p}{2} < 0$. 
Then \[
H'(z) = 2a \left(1 - \frac{p\gamma_p}{2}\right) z^{-\frac{p\gamma_p}{2}} + 4b \left(2 - \frac{p\gamma_p}{2}\right) z^{1 - \frac{p\gamma_p}{2}} < 0.
\]
Therefore,  $H$ is strictly decreasing on $(0, \infty)$ and
\[
\lim_{z \to 0^+} H(z) = +\infty, \quad \lim_{z \to \infty} H(z) = 0,
\]as shown in Figure~\ref{fig:HcFc}.  From this, for each $c > 0$, equation \eqref{eq:3} defines a unique $z(c) = H^{-1}(\gamma_p^2 p C_{3,p}^p c^{\beta}) > 0$,  which is continuous and strictly decreasing for $c\in(0,\infty)$.

\begin{figure}[H]
  \centering
  \includegraphics[width=0.38\textwidth]{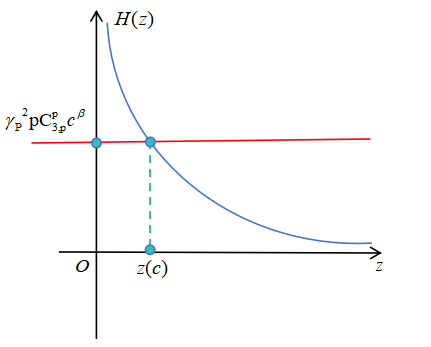}
    \includegraphics[width=0.38\textwidth]{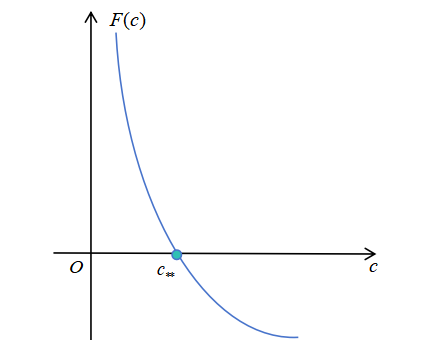}
  \caption{The figure of $H(c)$ and $F(c)$}
  \label{fig:HcFc}
\end{figure}

Define the function $F(c) := a(1 - \frac{2}{p\gamma_p}) z(c) + b \left(1 - \frac{4}{p\gamma_p}\right)   z(c)^2 - \tfrac{3}{2} c^2$. We analyze its behavior as $c \to 0^+$ and $c \to \infty$. As $c \to 0^+$, $c^{\beta} \to 0$ since $\beta >0$, so $z(c) \to \infty$ (since $H(z) \to 0$ as $z \to \infty$). Thus $F(c) \to +\infty$. As $c \to \infty$, $c^{\beta} \to \infty$, so $z(c) \to 0$ (since $H(z) \to \infty$ as $z \to 0^+$). Thus $F(c) \to -\infty$.
Since $F(c)$ is continuous  and strictly decreasing on $(0, \infty)$, it follows that there exists a unique $c_{**} > 0$ such that $F(c_{**}) = 0$ (see Figure~\ref{fig:HcFc}). Thus,  we have $F(c) >0$ if $c<c_{**}$.

\textbf{Case 2: $p=6$.}
When $p=6$, the function  
$\tilde{g}_c(r)$ 
 is independent of $c$. Consequently, $\tilde{g}_c$ has a unique global maximizer $$\widetilde{R} := \left(\frac{b + \sqrt{b^2 + 3 a \mathcal{S} ^{-3}}}{3 \mathcal{S} ^{-3}}\right)^{\frac{1}{2}}.$$
The corresponding maximum value is
$$\tilde{g}_c(\widetilde{R}) = \frac{2}{3}a\widetilde{R}^2+\frac{1}{3}b\widetilde{R}^4.$$
Define
$$\widetilde{c}_{**}:=\left(\frac{4}{9}a\widetilde{R}^2+\frac{2}{9}b\widetilde{R}^4\right)^{\frac{1}{2}}.$$
Then for $c<\widetilde{c}_{**}$, we have $\tilde{g}_c(\widetilde{R})> \frac{3 c^2}{2}$.

 \textbf{Conclusion.}
Finally, we define
\begin{equation*} 
c_{*}:=\begin{cases}
c_{**}, & 2<p<6,\\[4pt]
\widetilde{c}_{**}, & p=6,
\end{cases}
\qquad\text{and}\qquad
R_{*}:=\begin{cases}
R_{c}, & 2<p<6,\\[4pt]
\widetilde{R}, & p=6.
\end{cases}
\end{equation*}
For any $c \in (0,c_*)$,  we have $\max\limits _{{r>0}}\tilde{g}_c(r)=\tilde{g}_c(R_{*})>\frac{3}{2}c^{2}$.  Consequently, the equation $\tilde{g}_c(r) = \frac{3c^2}{2}$ admits exactly two positive solutions $R_0$ and $R_1$ satisfying $0 < R_0 < R_* < R_1$. Furthermore, the inequality $\tilde{g}_c(r) > \frac{3c^2}{2}$ holds precisely for $r \in (R_0, R_1)$. This completes the proof.
\end{proof}

\begin{lemma}\label{lem:c_02}
 For $ \frac{14}{3}<p\le6$, there exists    $c_0 > 0$ such that 
$R_*^2 >\frac{3}{2}c^2$ for all $c \in (0,c_0)$, where $R_*$ is defined in Lemma \ref{lem:critical_mass}.
\end{lemma}
\begin{proof}
For $\frac{14}{3}<p<6$,  from \eqref{eq:3}, we have  $R_{c}\to +\infty$ as $c\to 0^{+}$.  Hence, there exists  $c_{0,*} > 0$ such that 
$R_c^2 >\frac{3}{2}c^2$ for all $c \in (0,c_{0,*})$.
For $p=6$, set 
$$\widetilde{c}_{0,*}:=\left(\frac{2b + 2\sqrt{b^2 + 3 a \mathcal{S} ^{-3}}}{9 \mathcal{S} ^{-3}}\right)^{\frac{1}{2}}. $$ 
Thus, we have
$\widetilde{R}^2 >\frac{3}{2}c^2$ for all $c \in (0,\widetilde{c}_{0,*})$.
Then the conclusion of Lemma \ref{lem:c_02} follows by taking 
\begin{equation*} c_{0}:=\begin{cases}
c_{0,*}, & 2<p<6,\\[4pt]
\widetilde{c}_{0,*}, & p=6,
\end{cases}
\end{equation*}
  The proof  is complete.
\end{proof}

For $k>0$, we introduce the following set
\begin{align*}\label{gai-7}
   V_{k}(c)  := \{u\in S(c) \mid \|\nabla u\|_2 < k\},
\end{align*}
which serves as a natural constraint set for local minimization,
and consider the corresponding minimization problem
\begin{align*}
       m_{k}(c) := \inf_{u \in V_{k}(c)} I(u).
\end{align*}

 \begin{lemma}\label{lem:kongji}
Let $\frac{14}{3} < p \leq 6$.  
Then for $c\in(0,c_*)$, where $c_{*}$ is defined in Lemma~\ref{lem:critical_mass}, the following properties hold.

\begin{enumerate}
    \item[(i)] The Pohozaev manifold  $ \Lambda_{c}$ consists of two nonempty components $\Lambda_{c}^+$ and $ \Lambda_{c}^-$ satisfying
    \[
    \Lambda_{c}^+\subset V _{R_*}(c),    \qquad  \text{and}
    \qquad
     \Lambda_{c}^-\subset \{ u \in S(c) \mid \|\nabla u\|_2 > R_* \},
    \]
where $R_*$ is defined in Lemma~\ref{lem:critical_mass}. 

    \item[(ii)]     $ \Lambda_{c}^0  = \emptyset$.

    \item[(iii)]  $ \Lambda_{c}$ is a $C^1$ submanifold of $S(c)$ of codimension one.
\end{enumerate}
\end{lemma}

\begin{proof}
(i) Fix $u \in S(c)$, a direct computation yields \[
\Phi_u(s) = k(s)- \frac{3 c^2}{2}s + \frac{ c^2}{2}-\frac{1}{2}\int_{\mathbb{R}^3} u^2 \log u^2\,\mathrm{d}x,\ s\in \mathbb{R},
\]
where
\[
k(s) := \frac{a}{2}e^{2s}\|\nabla u\|_2^2 + \frac{b}{4}e^{4s}\|\nabla u\|_2^4 - \frac{1}{p}e^{p\gamma_p s}\|u\|_p^p,
\quad \gamma_p = \frac{3(p-2)}{2p} \in  (0,1].
\]
The first and second order derivatives of the fiber map are 
\[
\Phi_u'(s) = k'(s) - \frac{3 c^2}{2} \qquad  \text{and}\qquad \Phi_u''(s) = k''(s).
\]
Direct calculation shows that  function $\Phi_u(s)$  has at most two critical points.
Let $t = e^s$ and
\[
g(t) := a t^2\|\nabla u\|_2^2 + b t^4\|\nabla u\|_2^4 - \gamma_p t^{p\gamma_p}\|u\|_p^p.
\]
We observe that $k'(s) = g(t)$ and $k''(s) = tg'(t)$.

Set $r = t\|\nabla u\|_2 
= \|\nabla (s*u)\|_2> 0$. By the Gagliardo--Nirenberg inequality,
we obtain
$
g(t) \geq \tilde{g}_c(r),
$ where $\tilde{g}_c(r)$ is given in \eqref{Gai-2}.
Since $\Phi_u'(s) = g(e^s) - \frac{3  c^2}{2}$ and $g(e^s) \geq \tilde{g}(e^s\|\nabla u\|_2)$, 
we deduce from Lemma \ref{lem:critical_mass} that $\Phi_u'(s) > 0$ if $e^s\|\nabla u\|_2 \in (R_0, R_1)$. 
From the continuity of $\Phi_u'(s)$ with respect to $s$, $\lim\limits_{s\to-\infty}\Phi_u'(s) = -\frac{3  c^2}{2} < 0$ and $\lim\limits_{s\to+\infty}\Phi_u'(s) = -\infty$, we know that $\Phi_u'(s)$ must vanish at exactly two points $s_u^+ < s_u^-$, where $e^{s_u^+}\|\nabla u\|_2 < R_0<R_* <R_1<e^{s_u^-}\|\nabla u\|_2$.  Consequently, $s_u^+ *u \in  \Lambda_{c}^+ $ and $s_u^-  * u \in \Lambda_{c}^-$.
As every element of $\Lambda_{c}^\pm$ can be represented as $s_u^\pm * u$ for  $u \in S(c)$, the claimed inclusions follow.

For part (ii), suppose by contradiction that there exists a $u \in \Lambda_{c}^0$, i.e., $\Phi_u'(0) = \Phi_u''(0) = 0$. The condition $\Phi_u''(0) = 0$ implies $g'(1) = 0$, which means $t = 1$ is a critical point of $g$. By the Gagliardo--Nirenberg inequality, $\max \limits_{{t>0}}g(t)=g(1) \geq \max\limits _{{r>0}}\tilde{g}_c(r)=\tilde{g}_c(R_*)$.  However, for $c < c_*$ we have $\tilde{g}_c(R_*) > \frac{3c^2}{2}$, which contradicts $\Phi_u'(0) = g(1) - \frac{3 c^2}{2} = 0$. Therefore $\Lambda_{c}^0 = \emptyset$.

(iii) We now prove that $\Lambda_{c}$ is a $C^1$ submanifold of $S(c)$ of codimension 1. To this end, we define the smooth map
\[
L(u) := (Q(u), H(u)) : W \to \mathbb{R}^2,
\]
where $H(u) := \|u\|_2^2 - c^2$. Observe that
\[
\Lambda_{c}= \{ u \in W \mid L(u) = (0,0) \} = \{ u \in S(c) \mid Q(u) = 0 \}.
\]
Since $Q$ and $H$ are of $C^1$-class, the proof is complete provided we show that the differential
\[
(\mathrm{d}Q(u), \,\mathrm{d}H(u)) : W \to \mathbb{R}^2
\]
is surjective for every $u \in H^{-1}(0) \cap Q^{-1}(0)$.

Assume, for contradiction, that there exists a non-zero vector $(\eta_1, \eta_2) \in \mathbb{R}^2$ such that the rows of $DL(u)$ are linearly dependent, which means that there exists a non-zero constant $\eta \in \mathbb{R}$ satisfying
\begin{equation}\label{eq:linear-dependence}
\,\mathrm{d}Q(u)[\phi] = \eta \, \mathrm{d}H(u)[\phi] \quad \text{for all } \phi \in W.
\end{equation}
Define
\[
\psi(x) := \frac{\mathrm{d}}{\mathrm{d}s}\bigl(s * u(x)\bigr)\Big|_{s=0} = \frac{3}{2}u(x) + x \cdot \nabla u(x).
\]
Clearly $\psi \in W$, and a standard computation shows that $\psi \in T_u S(c)$, i.e., $\mathrm{d}H(u)[\psi] = 0$. Indeed, by integration by parts we obtain
\begin{align*}
\,\mathrm{d}H(u)[\psi] &= 2\int_{\mathbb{R}^3} u(x)\psi(x)\,\mathrm{d} x= 2\int_{\mathbb{R}^3} u(x)\left(\frac{3}{2}u(x) + x \cdot \nabla u(x)\right)\,\mathrm{d}x \\
&= 3\int_{\mathbb{R}^3} u^2(x)\,\mathrm{d} x+2\int_{\mathbb{R}^3} u(x)(x \cdot \nabla u(x))\,\mathrm{d} x\\
& = 0,
\end{align*}
where we have used  $\int_{\mathbb{R}^3} u(x)(x \cdot \nabla u(x))\,\mathrm{d}x = -\frac{3}{2}\int_{\mathbb{R}^3} u^2(x)\,\mathrm{d}x$.

Substituting $\phi = \psi$ into \eqref{eq:linear-dependence} and using $\mathrm{d}H(u)[\psi] = 0$, we deduce  $\mathrm{d}Q(u)[\psi] = 0$. From this and
\[
\mathrm{d}Q(u)[\psi] = \frac{\mathrm{d}}{\mathrm{d}s}\Big|_{s=0} Q(s * u) = \Phi_u''(0),
\]
we have $\Phi_u''(0) = 0$, which means $u \in \Lambda_{c}^0$. But this contradicts part (i) of the lemma, which guarantees that $\Lambda_{c}^0 = \emptyset$ for $0 < c < c_*$. Therefore, our assumption was false, and $DL(u)$ must be surjective for every $u \in \Lambda_{c}$. The proof is complete.
\end{proof}

With the structure of $\Lambda_{c}$ clarified, we now turn to a detailed analysis of the fiber map $\Phi_u(s)$.
\begin{lemma}\label{lem:critical-points}
Assume $\frac{14}{3} < p \leq 6$ and $c\in(0,c_*)$, where $c_{*}$ is defined in Lemma~\ref{lem:critical_mass}.
Then for each \( u \in S(c) \), the fibering map
\(
\Phi_u(s)
\)
has exactly two critical points \( s_u^+, s_u^- \in \mathbb{R} \) with \( s_u^+ < s_u^-\). Moreover,
\begin{enumerate}
\item[(i)] \( s_u^+ \) is the unique  local  minimizer for \( \Phi_u(s) \) with $\|\nabla  (s_u^+*u)\|_2<R_* $ (i.e.,  \(  s_u^+ * u \in \Lambda_{c}^+\)), and
\quad \( s_u^- \) is the unique  local  maximizer for \( \Phi_u(s) \) with $\|\nabla  (s_u^-*u)\|_2>R_* $ (i.e.,  \( s_u^- * u\in \Lambda_{c}^- \)), where $R_* $ is defined in Lemma~\ref{lem:critical_mass}.
\item[(ii)] The mappings $u \mapsto s_u^+$ and $u \mapsto s_u^-$ are of class $C^1$ on $S(c)$.
\end{enumerate}
\end{lemma}

\begin{proof}

The conclusion of (i) follows directly from Lemma \ref{lem:kongji}.
The $C^1$ regularity of the maps $u \mapsto s_u^\pm$ is obtained by applying the implicit function theorem to $\Phi_u'(s) = 0$, which is valid since $\Phi_u''(s_u^\pm) \neq 0$.
\end{proof}

Next, we show that  $m_{{c}}^\pm$ defined by \eqref{gai-6} can be taken over  the radial functions without loss of generality. Define
\begin{equation}\label{g-11}
 m_{r,c}^\pm := \inf_{u \in \Lambda_{r,c}^\pm} I(u)
\end{equation}
and
\[
\Lambda_{r,c}^\pm := \Lambda_c^\pm \cap H_r^1(\mathbb{R}^3).
\]

\begin{lemma}\label{lem:m+-}
Let $\frac{14}{3} < p \leq 6$ and $c\in(0,c_*)$, where $c_{*}$ is defined in Lemma~\ref{lem:critical_mass}. Then
\[
m_{{c}}^\pm = m_{r,c}^\pm.
\]

\end{lemma}

\begin{proof}
 Inequality $m_{r,c}^\pm  \geq m_{{c}}^\pm$ is trivial. For the reverse, given $u \in S(c)$, let $u^*$ be its Schwarz symmetric rearrangement. Using an argument similar to that in the proof of Lemma \ref{lem:m(c)=mr(c)}, we obtain  $I(u^*) \leq I(u)$.
By Lemma~\ref{lem:critical-points}, for each $ u \in S(c) $,
 one has $ s_u^+ \ast u \in \Lambda_c^+ $ and $ s_u^- \ast u \in \Lambda_c^- $ with $s_u^+<s_u^- $. To this end, we give the following equivalent characterizations of $m_{{c}}^\pm$ 
\begin{align}
\inf_{u \in \Lambda_c^+} I(u) = \inf_{u \in S(c)} \min_{-\infty < s \leq s_u^+} I(s*u), \ \text{and}\
\inf_{u \in \Lambda_c^-} I(u) = \inf_{u \in S(c)} \max_{s_u^+ < s \leq s_u^-} I(s*u).\label{lem:m+-1}
\end{align}
Recall that
\[
\lim_{s \to -\infty} \Phi_{u^*}'(s) \leq \lim_{s \to -\infty} \Phi_u'(s) = +\infty,
\]
\[
\Phi_{u^*}''(s) \leq \Phi_u''(s), \quad \forall s \in \mathbb{R},
\]
and
 \[
 0=\Phi_{u^*}'(s_{u^*}^{+}) \leq  \Phi_u'(s_{u^*}^{+}),
\]
which implies $-\infty < s_u^+ \leq s_{u^*}^+  \leq s_u^-$.
Similarly, it follows from
\[
 0=\Phi_{u^*}'(s_{u^*}^{-}) \leq  \Phi_u'(s_{u^*}^{-})
\]
that $-\infty < s_u^+ \leq s_{u^*}^- \leq s_u^-$.
Combining the above analysis, we get $-\infty < s_u^+ \leq s_{u^*}^+ < s_{u^*}^- \leq s_u^-$.
From this  and $$I((s*u)^*) \leq I(s*u)\ \ \text{for \ all} \ s \in \mathbb{R}, $$ we obtain
\begin{align*}
\min_{-\infty < s < s_{u^*}^+} I((s*u)^*) \leq \min_{-\infty < s < s_u^+} I(s*u), \ \text{and}\
\max_{s_{u^*}^+ < s \leq s_{u^*}^-} I((s*u)^*) \leq \max_{s_u^+ < s \leq s_u^-} I(s*u).
\end{align*}
Together with the above  and  \eqref{lem:m+-1}, we have  $$\inf_{u \in \Lambda_{r,c}^+} I(u) \leq \inf_{u \in \Lambda_c^+} I(u),\ \text{and}\ \inf_{u \in \Lambda_{r,c}^-} I(u) \leq \inf_{u \in \Lambda_c^-} I(u).$$The proof is complete.
\end{proof}

We consider the minimization problems
\begin{align}\label{gai-mRc}
     m _{R_*} (c) := \inf_{u \in V _{R_*}(c)} I(u),
\end{align}
and \begin{align}\label{mRcr}
    m _{r,R_*}(c) :=\inf\limits_{u\in V _{R_*}(c)\cup H^1_r(\mathbb{R}^3)}I(u),
\end{align}
where $R_*$  is defined in Lemma~\ref{lem:critical_mass}. 

As in the proof of Lemma \ref{lem:m(c)=mr(c)}, we have the following lemma.
\begin{lemma}\label{lem:mRcr}
    For $\frac{14}{3} < p \leq 6$ and $c\in(0,c_*)$, where $c_{*}$ is defined in Lemma~\ref{lem:critical_mass}, we have
\[
 m _{R_*}(c)=m _{r,R_*}(c).
\]
\end{lemma}

In what follows, we establish a lemma that connects the minimization problem on $\Lambda_c^+$ with the corresponding one on the bounded set $V _{R_*}(c)$. This equivalence allows us to work within a compact framework.

\begin{lemma}\label{lem:mAc=m+}
Assume  $\frac{14}{3} < p \leq 6$ and $c\in(0,c_*)$, where $c_{*}$ is defined in Lemma~\ref{lem:critical_mass}. Then the following properties hold.

\begin{enumerate}
\item[(i)]  $m^+_c = m _{R_*}(c).$
\item[(ii)]  
  For any $c\in(0,c_*) $, there exists a constant $\rho > 0$ 
small enough such that  
\[
m _{R_*}(c)< \inf_{\overline{V _{R_*}(c)} \setminus V _{R_*-\rho}(c)} I(u),
\]
where $\overline{V _{R_*}(c)}$ denotes the closure of $V _{R_*}(c)$ and
\[
\overline{V _{R_*}(c)}\setminus V _{R_*-\rho}(c) = \{u \in S(c) \mid R_*-\rho \le\|\nabla u\|_2 < R_* \}.
\]
\end{enumerate}
\end{lemma}
\begin{proof}
(i) By Lemmas \ref{lem:kongji} and  \ref{lem:critical-points}, for any $u \in \Lambda_{c}^+$ we have $s_u^+ * u \in \Lambda_{c}^+ \subset V _{R_*}(c)$, which immediately yields $m _{R_*}(c) \leq m^+_{c}$.

Conversely, for any $u \in V _{R_*}(c)$, the mapping $s \mapsto \|\nabla (s*u)\|_2^2 = e^{2s} \|\nabla u\|_2^2$ is continuous and strictly increasing. Thus, there exists a unique $s_u^1 \in \mathbb{R}$ such that $\|\nabla (s_u^1*u)\|_2^2 = R_*$. Then, one has
\[
\|\nabla (s_u^+*u)\|_2^2 < \|\nabla (s_u^1*u)\|_2^2 = R_* < \|\nabla (s_u^-*u)\|_2^2,
\]
which implies $s_u^+ < s_u^1 < s_u ^-$. Consequently, we have
\begin{align*}
I(s_u^+*u) &= \min \left\{ I(s*u) \mid s \in \mathbb{R} \text{ and } s < s_u^1 \right\} \\&= \min \left\{ I(s*u) \mid s \in \mathbb{R} \text{ and } \|\nabla (s*u)\|_2^2 < R_* \right\} \leq I(u).
\end{align*}
 It follows from  the above that
\[
m^+_c \leq m _{R_*}(c),
\]
which establishes the first part of the lemma.

(ii) 
Fix $c>0$ and any $u\in \Lambda_{c}^+$. By   Lemma  \ref{lem:kongji},   one has    $\|\nabla   u \|_2<R_{0}<R_* $.
Choose $\rho>0$ such that $\rho<R_{*}-R_{0}$. Then,
one has  $\|\nabla   u\|_2<R_{0}<R_{*}-\rho<R_* $. Hence, 
  $ \Lambda_{c}^+ \subset V _{R_*-\rho}(c)$. This completes the proof.
\end{proof}

\section{$L^{2}$-supercritical case: $\frac{14}{3} < p < 6$}
\setcounter{equation}{0}

We are now prepared to prove the existence of  solutions for $\frac{14}{3} < p < 6$. The following lemma establishes the existence of a ground state  as a local minimizer.

\begin{lemma}\label{lem:mAc}
Assume $\frac{14}{3} < p < 6$ and $c\in(0,c_{**})$, where $c_{**}$ is defined in Lemma~\ref{lem:critical_mass}. Then  $m_{R_c}(c)$ defined in    \eqref{gai-mRc} with $R_{*}=R_c$ 
is achieved by a  positive and  radially symmetric function $u_c^+ \in V _{R_c}(c)$. Furthermore,
$u_{c}^{+}$ is a ground state for $I|_{S_c}$.
\end{lemma}

\begin{proof}
From Lemma \ref{lem:mRcr},  it remains only to  prove that $m _{r,R_c}(c)$ defined by \eqref{mRcr} with $R_{*}=R_{c}$ can be achieved by a function in $ V _{R_c}(c)\cap H^1_r(\mathbb{R}^3)$, which is a positive and radially symmetric
function.
Let $\{u_n\} \subset V _{R_c}(c)\cap H^1_r(\mathbb{R}^3)$ be a minimizing sequence for $m _{r,R_c}(c)$. It is easily seen that $\{u_n\} $ is bounded
in $H_r^1(\mathbb{R}^3)$. 
By Lemma~\ref{lem:mAc=m+} and an argument similar to that in the proof of 
 Theorem \ref{thm:main1}, one has 
$u_n  \to  u_c^{+} $ strongly in $W_{r}$ as $n \to \infty$ and  $u_c^{+}\in V _{R_c}(c)$. Hence, $m_{R_c}(c)$ 
is achieved by $u_c^{+}$.

Using arguments similar to those in Theorem \ref{thm:main1},     there exists a  $\lambda_c^+\in\mathbb{R} $ such that $(\lambda_c^+,u_c^+)$ solves  problem \eqref{eq:main-1}--\eqref{eq:main-2}, where  $u_c^+$ is positive and radially symmetric. 
Since any critical point of $I|_{S(c)}$ lies in $\Lambda_{c}$ and $m _{R_c}(c)=m^+_{c} \leq m^-_{c}$, we see that $ u_c^{+}$ is a ground state for $I|_{S(c)}$. The proof is complete.
\end{proof}

Finally, we establish the existence of the second solution.

\begin{lemma}\label{lem:m+-achieved}
Assume $\frac{14}{3} < p < 6$ and $c\in(0,c_{**})$, where $c_{**}$ is defined in Lemma~\ref{lem:critical_mass}.
Then  $m^-_{c}$ defined in \eqref{gai-6} is achieved by a   positive and radially symmetric function $u_c^- \in S(c) \cap H_r^1(\mathbb{R}^3)$.
\end{lemma}

\begin{proof}
In view of Lemma~\ref{lem:m+-}, it remains to prove the existence of a positive and  radially symmetric minimizer for $m^-_{r,c}$ in $ \Lambda^-_{r,c}$.
Let $\{u_n\} \subset \Lambda^-_{r,c}$ be a minimizing sequence of
$m^-_{r,c}$ defined in \eqref{g-11}. Hence  $s_{u_{n}}^-=0$.
By the Gagliardo-Nirenberg  inequality,   property (ii) of $B(u)$ and $Q(u_n)=0$, we obtain
\begin{align*}
&\ \ \ \ m^-_{r,c}+1+o_{n}(1)\|  u_{n}\|_W \ge I(u_{n})-\frac{1}{p}Q(u_{n}) \\ 
&\geq \left( \frac{1}{2} - \frac{1}{p \gamma_p} \right) a\|\nabla u_{n}\|_2^2 + \left( \frac{1}{4} - \frac{1}{p \gamma_p} \right) b\|\nabla u_{n}\|_2^4 + \frac{1}{2} \int_{\mathbb{R}^3} \left( A(|u_{n}|) - B(|u_{n}|) \right) \mathrm{d}x \nonumber
\\
&\geq \left( \frac{1}{2} - \frac{1}{p \gamma_p} \right) a\|\nabla u_{n}\|_2^2 + \left( \frac{1}{4} - \frac{1}{p \gamma_p} \right) b\|\nabla u_{n}\|_2^4 -\frac{K_q}{2} \|u_n\|_q^q \nonumber
\\
&\geq  \left( \frac{1}{4} - \frac{1}{p \gamma_p} \right) b\|\nabla u_{n}\|_2^4   -\frac{ K_q C^{q}_{3, q}}{2} \|\nabla u_{n}\|_2^{q\gamma_q} c^{ (1 - \gamma_q)q}.
\end{align*}
Thus, $\{u_n\}$ is bounded in $H_{r}^1(\mathbb{R}^3)$  since $q\gamma_q<2$,  and $\int_{\mathbb{R}^3} A(|u_n|) \mathrm{d}x$ is bounded. By Lemma \ref{lem:Orlicz-V}, $\{u_n\}$ is bounded in $W_r$. Thus, up to a subsequence, there exists a $u_c^{-} \in W_r$ such that $u_n \rightharpoonup u_c^{-} $ weakly in $W_r$. By Lemma \ref{lem:compact-embedding}, we have $u_n \to u_c^{-}$ strongly in $L^2(\mathbb{R}^3)$ and a.e. in $\mathbb{R}^3$.
Consequently,
we know that $u_n \to u$ strongly in $L^l(\mathbb{R}^3)$ for all $l \in \left(2, 6\right)$ by Lemma \ref{lem:G-N}, 
which in turn implies that
$B(|u_n|) \to B(|u|) $ strongly in $L^1(\mathbb{R}^3)$.
Combining this with  the weak lower semi-continuity, we have $Q(u_c^{-}) \leq  0$ and $0=\Phi_{u_c^{-} }'(s_{u_c^{-}}^-) \leq \liminf\limits_{n \to \infty}\Phi_{u_n}'(s_{u_c^{-}}^-)$. Hence, $s_{u_n}^+\leq s_{u_c^{-}}^-\leq s_{u_n}^-=0$. Recalling from Lemma~\ref{lem:critical-points} that $s_{u_n^{-}}=0$ is the unique local maximizer of $\Phi_{u_n}(s)$ and using the weak lower semi-continuity again, one has 
\[
m^-_{r,c} = I(u_n)+ o_n(1)  \geq I(s_{u_c^{-}}^-*u_{n}) + o_n(1) \geq I(s_{u_c^{-}}^-*u_{c}^{-})  + o_n(1) \geq m^-_{r,c},\  \text{as}\  n \to \infty.
\]
Therefore, $s_{u_c^{-}}^-=0$, $I(u_{c}^-) = m^-_{r,c}$ and $A(|u_n|) \to A(|u_{c}^-|)$ strongly in $L^1(\mathbb{R}^3)$. By Lemma \ref{lem:Orlicz-V},   $u_n \to u_{c}^-$ in $V$ as $n \to \infty$, and hence in $W_r$,  which implie that $m^-_{r,c}$ is attained by $u_c^-$.

Similar to the proof of  Theorem \ref{thm:main1}, 
 there exists a $\lambda _c^{-}\in \mathbb{R}$ such that $(\lambda_c^+,u_c^+)$ solves  problem \eqref{eq:main-1}--\eqref{eq:main-2}, where  $u_c^-$ is positive and radially symmetric.
 The proof is complete.
\end{proof}

\begin{proof}[Proof of Theorem \ref{thm:main2}]
The case for $m^+_{c}$ follows from Lemmas \ref{lem:mAc=m+} and \ref{lem:mAc}, while the case for $m^-_{c}$ follows from Lemma \ref{lem:m+-achieved}. \end{proof}

\section{The Sobolev critical case: \(p=6\)}
\setcounter{equation}{0}
In this section, we address the Sobolev critical case $p=6$, where the Sobolev embedding $H^{1}(\mathbb{R}^{3})\hookrightarrow L^{6}(\mathbb{R}^{3})$ fails to be compact. 

\subsection{Existence of a ground state}

We establish a local compactness condition for the energy functional associated with problem \eqref{eq:main-1}--\eqref{eq:main-2}, which will be essential for constructing  the ground state solution.

\begin{lemma}\label{prop:PS-critical}
Let $\{u_n\} \subset S(c) \cap H_r^1(\mathbb{R}^3)$ be a Palais--Smale sequence for $I|_{S(c)}$ at level $\widetilde{m}$ satisfying $Q(u_n) \to 0$. Then, up to a subsequence, one of the following alternatives holds:

\begin{itemize}
\item[(i)] $u_n \rightharpoonup u$ weakly but not strongly in $W_r$, where $u$ satisfies the limit equation
\[
-(a + Kb)\Delta u = \lambda u + |u|^4u + u\log u^2 \quad \text{in } \mathbb{R}^3
\]
for some $\lambda \in \mathbb{R}$ and $K := \lim\limits_{n\to\infty}\|\nabla u_n\|_2^2>0$. Moreover,
\[
\widetilde{m} \geq I(u) + \frac{a + b\|\nabla u\|_2^2}{3}M + \frac{b}{12}M^{2},
\]
with $M$ defined by
\begin{equation}\label{M}
M:=\frac{\mathcal{S}}{2}\Big(b\mathcal{S}^2 + \sqrt{  b^2\mathcal{S}^4+4 (a+b\|\nabla u\|_2^2)\mathcal{S}}\Big).
\end{equation}

\item[(ii)] $u_n \to u$ strongly in $W_r$ for some $u \in W_r$. Consequently, $I(u) = \widetilde{m}$ and $(\lambda, u)$ solves problem \eqref{eq:main-1}--\eqref{eq:main-2} for some $\lambda \in \mathbb{R}$.
\end{itemize}
\end{lemma}

\begin{proof}

We divide the proof into four steps.

\textbf{Step 1: Boundedness and weak convergence.}
Similar to the proof of Lemma \ref{lem:m+-achieved}, one has  $\{u_n\}$ is bounded in $W_r$. Therefore, up to a subsequence, one has
\[
u_n \rightharpoonup u \quad \text{weakly in } W_r, \qquad
u_n \to u \quad \text{strongly in } L^2(\mathbb{R}^3), \  \text{and}\ \
u_n \to u \quad \text{a.e. in } \mathbb{R}^3.
\]
Consequently,
  $u_n \to u$ strongly in $L^l(\mathbb{R}^3)$ for all $l \in \left(2, 6\right)$ by Lemma \ref{lem:G-N}. 

\textbf{Step 2: Existence of Lagrange multipliers.}
From $I'|_{S(c)}(u_n) \to 0$ and the Lagrange multiplier rule (see \cite[Lemma 3]{Berestycki1983}), there exists a sequence $\{\lambda_n\} \subset \mathbb{R}$ such that for every $
 \phi \in W$,
\begin{align}\label{eq:lagrange2}
&(a + b\|\nabla u_n\|_2^2)\int_{\mathbb{R}^3} \nabla u_n  \nabla\phi\,\mathrm{d}x - \int_{\mathbb{R}^3} \lambda_n u_n \phi\,\mathrm{d}x \nonumber\\
 = &\int_{\mathbb{R}^3} |u_n|^4 u_n \phi\,\mathrm{d}x+ \int_{\mathbb{R}^3} u_n \phi \log u_n^2\,\mathrm{d}x + o_n(1)\|\phi\|_{W}.
\end{align}
Taking $\phi = u_n$ and using the boundedness of $\{u_n\}$ in $W_r$, we see that $\{\lambda_n\}$ is bounded. Thus, up to a subsequence, $\lambda_n \to \lambda$ for some $\lambda \in \mathbb{R}$.

\textbf{Step 3: Energy decomposition.}
Since  $Q(u_n) \to 0$,  we obtain
\[
a\|\nabla u_n\|_2^2 + b\|\nabla u_n\|_2^4 = \|u_n\|_6^6 + \frac{3  c^2}{2} + o_n(1) \geq \frac{3  c^2}{2} + o_n(1).
\]
This shows that $\{\|\nabla u_n\|_2\}$ is  away from zero.  Thus \(K> 0\). Passing to the limit in \eqref{eq:lagrange2} yields
\begin{equation}\label{Lem:pro-eq2}
(a+Kb)\int_{\mathbb{R}^{3}}\nabla u\nabla\phi\,\mathrm{d}x-\lambda\int_{\mathbb{R}^{3}}u\phi\,\mathrm{d}x
-\int_{\mathbb{R}^{3}}|u|^{4}u\phi\,\mathrm{d}x
-\int_{\mathbb{R}^{3}}u\phi\log u^{2}\,\mathrm{d}x
=0,
\qquad\forall\phi\in W.
\end{equation}
Hence \(u\) satisfies \(-(a+Kb)\Delta u = \lambda u + |u|^4u +  u\log u^2\), and the associated Pohozaev identity reads
\begin{equation}\label{Lem:pro-eq3}
\widetilde{Q}(u) := (a+Kb)\|\nabla u\|_2^2 - \|u\|_6^6 - \frac{3 c^2}{2} = 0.
\end{equation}
Setting \(v_n := u_n - u\). Then, we have \(v_n \rightharpoonup 0\) in \(W_{r}\),  \(v_n \to 0\) in \(L^l(\mathbb{R}^3)\) for $l\in(2,6)$ and
\begin{equation}\label{Lem:pro-eq3a11}
\|\nabla u_n\|_2^2 = \|\nabla u\|_2^2 + \|\nabla v_n\|_2^2 + o_n(1),\ \text{as}\  n\to+\infty.
\end{equation}
It follows from the property (ii) of $B(u)$, Lemma \ref{lem;duishu}  and  the Brézis-Lieb lemma that
\begin{align}
 \int_{\mathbb{R}^3} u_n^2 \log u_n^2\,\mathrm{d}x
  = &\int_{\mathbb{R}^3} v_n^2 \log v_n^2\,\mathrm{d}x+\int_{\mathbb{R}^3} u^2 \log u^2 \,\mathrm{d}x+o_{n}(1)\nonumber\\
 \le& K_q\|v_n\|_q^q+\int_{\mathbb{R}^3} u^2 \log u^2 \,\mathrm{d}x+o_{n}(1)\nonumber\\
 =&\int_{\mathbb{R}^3} u^2 \log u^2 \,\mathrm{d}x+o_{n}(1)\label{eq1-lem;duishu}
\end{align} and
\begin{equation}\label{Lem:pro-eq3a1}
\|u_n\|_6^6 = \|u\|_6^6 + \|v_n\|_6^6 + o_n(1),\ \text{as}\  n\to+\infty.
\end{equation}
The condition \(Q(u_n) = o_n(1)\) can be rewritten as
\begin{equation}\label{Lem:pro-eq3--1}
Q(u_n) = (a+Kb)\|\nabla u_n\|_2^2 - \|u_n\|_6^6 - \frac{3 c^2}{2} =o_n(1),\ \text{as}\  n\to+\infty.
\end{equation}
From  \eqref{Lem:pro-eq3}, \eqref{Lem:pro-eq3a11}, \eqref{Lem:pro-eq3a1}, and \eqref{Lem:pro-eq3--1}, we obtain
\begin{equation}\label{Lem:pro-eq41}
\|v_n\|_6^6 = (a+Kb)\|\nabla v_n\|_2^2 + o_n(1),\ \text{as}\  n\to+\infty.
\end{equation}
Consequently,
\begin{equation}\label{Lem:pro-eq41-1}
\lim_{n\to\infty}\|v_n\|_6^6 \geq \lim_{n\to\infty}[(a+b\|\nabla u\|_2^2)\|\nabla v_n\|_2^2 + b\|\nabla v_n\|_2^4].
\end{equation}
Let \(\ell := \lim\limits_{n\to\infty}\|\nabla v_n\|_2^2\). The Sobolev inequality \(\mathcal{S}\|u\|_6^2 \leq \|\nabla u\|_2^2\) yields
\[
\mathcal{S}^{-3}\ell^{3} \geq (a+b\|\nabla u\|_2^2)\ell + b\ell^{2}.
\]
This implies either \(\ell = 0\) or
$\ell\geq M$, where $M$ is defined by \eqref{M}.

\textbf{Step 4:  Compactness results for  the Palais-Smale sequences.} Two possibilities emerge.
\begin{enumerate}
\item[(i)]   If \(\ell \geq M\), then using \eqref{Lem:pro-eq3a11},\eqref{eq1-lem;duishu},\eqref{Lem:pro-eq3a1}, and
\eqref{Lem:pro-eq41}, we obtain
\begin{align*}
\widetilde{m} &= \lim_{n\to+\infty}I(u_n)\\
&= I(u) + \lim_{n\to+\infty}\left[\frac{a+b\|\nabla u\|_2^2}{2}\|\nabla v_n\|_2^2 + \frac{b}{4}\|\nabla v_n\|_2^4 - \frac{1}{6}\|v_n\|_6^6-\frac{1}{2}\int_{\mathbb{R}^3} v_n^2 \log v_n^2\,\mathrm{d}x\right]\\
&= I(u) + \lim_{n\to+\infty}\left[\frac{a+b\|\nabla u\|_2^2}{3}\|\nabla v_n\|_2^2 + \frac{b}{12}\|\nabla v_n\|_2^4-\frac{1}{2}\int_{\mathbb{R}^3} v_n^2 \log v_n^2\,\mathrm{d}x\right]\\
&\geq I(u) + \frac{a+b\|\nabla u\|_2^2}{3}M+ \frac{b}{12}M^2.\end{align*}
This corresponds to alternative (i).
\item[(ii)] If \(\ell = 0\), then \(u_n \to u\) strongly in \(W_{r}\) and \(L^6(\mathbb{R}^3)\), yielding alternative (ii).
\end{enumerate}
The proof is complete.
\end{proof}

Employing the compactness framework established above, we consider the existence of a ground state solution.

\begin{proof}[Proof of Theorem \ref{thm:main3}]

Let \(\{v_n\} \subset V _{\widetilde{R}}(c)\) be a minimizing sequence for \(m _{\widetilde{R}}(c)\)  defined by   \eqref{gai-mRc} with $R_{*}=\widetilde{R}$.
 By Lemma \ref{lem:mRcr}, we   assume \(\{v_n\} \subset H_r^1(\mathbb{R}^3) \cap V _{\widetilde{R}}(c)\) is radial for each \(n\).
Lemmas \ref{lem:critical-points} and \ref{lem:mAc=m+} guarantee the existence of \(s_{v_n}^+\) such that \(s_{v_n}^+ * v_n \in \Lambda_{r,c}^+\) with \(\|\nabla(s_{v_n}^+ * v_n)\|_2 < \widetilde{R} - \rho\), and
\[
I(s_{v_n}^+ * v_n) = \min\left\{I(s * v_n) : s \in \mathbb{R},\ \|\nabla(s * v_n)\|_2 < \widetilde{R}\right\} \leq I(v_n),
\]
where $\rho>0$ is small enough.
Thus, \(\{w_n := s_{v_n}^+ * v_n\}\) constitutes a minimizing sequence for \(m _{\widetilde{R}}(c)\)  and \(Q(w_n) = 0\).

Applying  Ekeland's variational principle, we obtain another minimizing sequence \(\{u_n\} \subset S(c) \cap H_r^1(\mathbb{R}^3)\) for \(m _{\widetilde{R}}(c)\) satisfying
\begin{equation}\label{eq-gai-6.10}
I(u_n) \leq I(w_n), \quad \|u_n - w_n\|_{ W_{r}} \leq \frac{\rho}{2},
\end{equation}
and
\begin{equation}\label{Eke}
I(u_n) < I(w) + \frac{2}{n\rho}\|w - u_n\|_{ W_{r}}, \quad \forall \ w \in  S(c)\cap H_r^1(\mathbb{R}^3),\quad w \neq u_n.
\end{equation}
It follows that \(\{u_n\} \subset V _{\widetilde{R}}(c)\). For any \(\widetilde{h} \in T_{u_n}S_r(c)\) with \(\|\widetilde{h}\|_{ W_{r}} = 1\), set
\[
z(t) = \frac{u_n + t\widetilde{h}}{\|u_n + t\widetilde{h}\|_2}c, \quad t \in \left(0, \frac{\rho}{2}\right].
\]
Then \(z(t) \in S(c)\) and \(\|z(t)\|_{W_{r}} < {\widetilde{R}}\). Substituting \(w = z(t)\) in \eqref{Eke} and letting \(t \to 0^+\) yields
\[
\langle I'(u_n), \widetilde{h} \rangle \geq -\frac{2}{n\rho}.
\]
Replacing \(\widetilde{h}\) by \(-\widetilde{h}\) gives \(\langle I'(u_n), \widetilde{h} \rangle \leq \frac{2}{n\rho}\). Thus, \(\{u_n\}\) is a Palais--Smale sequence for \(I|_{S(c)}\) at level \(m _{\widetilde{R}}(c)\), i.e., \(I(u_n) \to m _{\widetilde{R}}(c)\) and \(I'|_{S(c)}(u_n) \to 0\) $ \text{as}\  n\to+\infty$. From \eqref{eq-gai-6.10} and $\rho>0$ is small enough, we have  \(\|\nabla u_n\|_2 \leq \widetilde{R} - \frac{\rho}{2}<\widetilde{R}\) and \(Q(u_n) \to 0\) $\text{as}\  n\to+\infty$.   Therefore, \(\{u_n\} \) satisfies the assumptions of Lemma \ref{prop:PS-critical}.

We claim that (ii) in Lemma \ref{prop:PS-critical} occurs. If not, by an argument similar to that in the proof of Lemma \ref{prop:PS-critical}, we obtain  \(u_n \rightharpoonup  u_c^{+} \) weakly in \(W_{r}\) but not strongly, where \( u_c^{+}\)   satisfies
\[
m _{\widetilde{R}} (c)\geq I( u_c^{+}) + \frac{a + b\|\nabla  u_c^{+}\|_2^2}{3}M+ \frac{b}{12}M^2 > I( u_c^{+}).
\]
By the weak lower semi-continuity, we have   \( u_c^{+} \in V _{\widetilde{R}}(c)\). Consequently, \(m _{\widetilde{R}}(c) \leq I( u_c^{+})\), which yields a contradiction.
Therefore, \(u_n \to u_c^{+}\) strongly in \(W_{r}\), which implies that   $m _{\widetilde{R}}(c)$ is attained by $u_c^+$.

Using arguments similar to those in Theorem \ref{thm:main1} and Lemma \ref{lem:mAc},     there exists  a $\lambda_c^+\in\mathbb{R} $ such that $(\lambda_c^+,u_c^+)$ solves  problem \eqref{eq:main-1}--\eqref{eq:main-2}, where  $u_c^+$ is positive and  radially symmetric  ground state for $I|_{S_c}$. 
 The proof is complete.
\end{proof}

\subsection{Existence of a mountain pass sulution}
To obtain a second solution, we introduce an auxiliary functional and derive refined energy estimates on the Aubin-Talenti bubble to establish the compactness threshold.

\subsubsection{Auxiliary functional and its minimizer }

Following the ideas of \cite{ChenTang2025}, we introduce  a proper approach based on an auxiliary functional
 \begin{align}
\mathcal{J}(u) := &\frac{(b^2 \mathcal{S}^4 + 4a\mathcal{S})^{3/2}}{24} \left[\left(1 + \frac{4b}{b^2 \mathcal{S}^3 + 4a} \|\nabla u\|_2^2\right)^{3/2} - 1\right] \nonumber
\\&+ \left(\frac{a}{2} + \frac{b^2 \mathcal{S}^3}{4}\right) \|\nabla u\|_2^2 + \frac{b}{4} \|\nabla u\|_2^4 - \frac{1}{6} \|u\|_{6}^6 +\frac{1}{2}c^2- \frac{1}{2} \int_{\mathbb{R}^3} u^2 \log u^2 \,\mathrm{d}x,\label{gai-J}
\end{align}
 which will used to  control the mountain pass level of the energy functional of problem \eqref{eq:main-1}--\eqref{eq:main-2} to satisfy the local compactness condition.

The key observation is that for any $u\in S(c)$, we have
$\mathcal{J}(u)\ge I(u)$.
The associated Pohozaev functional for $\mathcal{J}$ is
\begin{align*}
P(u):=
{\Phi^{\mathcal{J}}_u}'(0)&=\Bigl[a+b\|\nabla u\|_2^2 +\frac{b^2\mathcal{S}^{3}}{2}+\frac{b\mathcal{S}}{2}\sqrt{b^2\mathcal{S}^4+4(a+b\|\nabla u\|_2^2 )\mathcal{S}}\Bigr]\|\nabla u\|_2^2 - \|u\|_6^6 - \frac{3 c^2}{2},
\end{align*}
where $\Phi^{\mathcal{J}}_u(s)=:\mathcal{J}(s*u)$.
Hence  any critical point of $\mathcal{J}|_{S(c)}$ satisfies $P(u)=0$.
Set
 $$\widetilde{\Lambda}_{c}^+:= \{ {\Phi^{\mathcal{J}}_u}'(0) = 0, \ \text{and}\ \ {\Phi^{\mathcal{J}}_u}''(0) > 0\}$$
and
\begin{align}\label{gai-8}
 \hat{m}_{{\widetilde{R}}}(c)=\inf\limits_{V  _{{\widetilde{R}}}(c)}\mathcal{J},\end{align}
where ${\widetilde{R}}>0$ is defined in Lemma \ref{lem:critical_mass}.

The following lemma give some basic properties for $\mathcal{J}(u)$.

\begin{lemma}\label{222}
Assume $p=6$ and $c\in(0,\widetilde{c}_{**})$, where $\widetilde{c}_{**}$ is defined in Lemma~\ref{lem:critical_mass}. Then we have

\begin{enumerate}
\item[(i)] $\widetilde{\Lambda}_{c}^+\subset V_{{\widetilde{R}}}(c)$.
 \item[(ii)]  $
\hat{m} _{{\widetilde{R}}}(c)< \inf\limits_{\overline{V _{{\widetilde{R}}}(c) }\setminus V _{{\widetilde{R}}-\rho}(c)} \mathcal{J}(u).
$
\end{enumerate}
\end{lemma}
\begin{proof}  For  any $c<\widetilde{c}_{**}$, by the Gagliardo-Nirenberg inequality, one has
\begin{align*}
P(s*u)= {\Phi^{\mathcal{J}}_u}'(s) 
 \ge a e^{2s} \|\nabla u\|_2^2
+ b e^{4s} \|\nabla u\|_2^4
- e^{6s}  \|u\|_{6}^6 - \frac{3c^{2}}{2}
 \ge k'(e^s\|\nabla u\|_2)- \frac{3c^{2}}{2},
\end{align*}
where the function $k$ is defined in Lemma \ref{lem:critical_mass}.
By arguments similar to those in Lemmas \ref{lem:critical_mass} and \ref{lem:kongji}, the geometric structure of  $\Phi^{\mathcal{J}}_u(s)$ is similar to that of  $\Phi_u(s)$.
Hence the conclusion is immediate.
\end{proof}

To obtain the existence of  solutions of $\hat{m} _{{\widetilde{R}}}(c)$, we introduce the following technical  condition
\begin{equation}\label{eq:technical}
\frac{4a^2}{b^3\mathcal{S}^3}+\frac{3a}{b}+\frac{3b\mathcal{S}^3}{4}\geq {\widetilde{R}}^2,
\end{equation}
where $\widetilde{R} $ is defined in Lemma \ref{lem:critical_mass}.

 \begin{lemma}\label{Lem2:solution}
Let $c\in(0,\widetilde{c}_{**})$ and assume the technical condition \eqref{eq:technical} holds, where $\widetilde{c}_{**}$ is defined in Lemma~\ref{lem:critical_mass}. 
Then $\hat{m} _{{\widetilde{R}}}(c)$ defined in \eqref{gai-8} is reached by a positive, radially symmetric and  non-increasing  function $\hat{u}_{c} \in  V _{\widetilde{R}}(c)$ with some  $\hat{\lambda}_{c} \in \mathbb{R}$.

\end{lemma}

\begin{proof}
Let $\{v_n\} \subset V _{{\widetilde{R}}}(c)$ be a minimizing sequence for $\hat{m} _{{\widetilde{R}}}(c)$. 
  Taking into account that
the Schwarz rearrangement of \(|v_n|\), up to a subsequence, we  get a new minimizing sequence
\(\{v_n\} \) which is nonnegative, radially symmetric and non-increasing in $r = |x|$. 
According to Lemma \ref{222} and the proof of 
 Theorem~\ref{thm:main3}, there exists    a  Palais--Smale sequence \(\{u_n\}\subset S(c) \cap H_r^1(\mathbb{R}^3) \) for \(\mathcal{J}|_{S(c)}\) at level \(\hat{m} _{{\widetilde{R}}}(c)\), i.e., \(\mathcal{J}(u_n) \to \hat{m} _{{\widetilde{R}}}(c)\) and \(\mathcal{J}'|_{S(c)}(u_n) \to 0\) and \(P(u_n) \to 0\). Moreover,  for each $n$, one has $u_n$ is   non-increasing in $r = |x|$ and 
\(\|\nabla u_n\|_2 \leq {\widetilde{R}} - \rho<{\widetilde{R}}\).

 \textbf{Step 1: Boundedness and weak convergence.} 
From \(P(u_n) \to 0\), the Gagliardo-Nirenberg inequality and
 property (ii) of $B(u)$, 
 we have
\begin{align*}
&\ \ \ \ \hat{m} _{{\widetilde{R}}}(c)+1+o_{n}(1)\|  u_{n}\|_W \ge \mathcal{J}(u_{n})-\frac{1}{6}P(u_{n})\\
&\ge \Bigl[\frac{a}{3}+\frac{b^2\mathcal{S}^3}{6}-\frac{b\mathcal{S}}{12} \sqrt{b^2\mathcal{S}^4+4(a+b\|\nabla u_{n}\|_2^2)\mathcal{S}}\Bigr]\|\nabla u_{n}\|_2^2\\
&
\quad+\frac{b}{12}\|\nabla u_{n}\|_2^4+\frac{1}{2}c^2-\frac{ K_q C^{q}_{3, q} }{2} \|\nabla u_{n}\|_2^{q\gamma_q} c^{ (1 - \gamma_q)q},
\end{align*}
which yields that $\{u_n\}$ is bounded in $W_r$. Therefore, up to a subsequence,
\[
u_n \rightharpoonup u \quad \text{weakly in } W_r, \qquad
u_n \to u \quad \text{strongly in } L^2(\mathbb{R}^3), \qquad
u_n \to u \quad \text{a.e. in } \mathbb{R}^3, 
\]Consequently,
we know that $u_n \to u$ strongly in $L^l(\mathbb{R}^3)$ for all $l \in \left(2, 6\right)$ by Lemma \ref{lem:G-N}. 
  By the pointwise convergence, one has  the limit function $u$  is radially symmetric and non-increasing with respect to $r = |x|$.

\textbf{Step 2: Existence of Lagrange multipliers.}
Since $\mathcal{J}'|_{S(c)}(u_n) \to 0$, by the Lagrange multiplier rule (see \cite[Lemma 3]{Berestycki1983}) there exists a sequence $\{\lambda_n\} \subset \mathbb{R}$ such that for every $\phi \in W$,
\begin{align}
&\Bigl[a+b\|\nabla u_n\|_2^2+\frac{b^2\mathcal{S}^{3}}{2}+\frac{b\mathcal{S}}{2}\sqrt{b^2\mathcal{S}^4+4(a+b\|\nabla u_n\|_2^2)\mathcal{S}}\Bigr]\int_{\mathbb{R}^3} \nabla u_n   \nabla\phi\,\mathrm{d}x \nonumber\\
 =  &\lambda_n  \int_{\mathbb{R}^3} u_n \phi\,\mathrm{d}x+\int_{\mathbb{R}^3} |u_n|^4 u_n \phi\,\mathrm{d}x + \int_{\mathbb{R}^3} u_n \phi \log u_n^2\,\mathrm{d}x + o_n(1)\|\phi\|_{W},\ \text{as}\  n\to+\infty .\label{eq1:Lem2:solution-1}
\end{align}
Taking $\phi = u_n$ and using the boundedness of $\{u_n\}$ in $W_r$, we see that $\{\lambda_n\}$ is bounded. Thus, up to a subsequence, $\lambda_n \to \lambda$ for some $\lambda\in \mathbb{R}$.

\textbf{Step 3: Energy decomposition.}
Since  $P(u_n) \to 0$,  we obtain
\begin{align*}
\Bigl[a+b\|\nabla u_n\|_2^2+\frac{b^2\mathcal{S}^{3}}{2}+\frac{b\mathcal{S}}{2}\sqrt{b^2\mathcal{S}^4+4(a+b\|\nabla u_n\|_2^2)\mathcal{S}}\Bigr]\|\nabla u_n\|_2^2
\geq  \frac{3c^{2}}{2}+ o_n(1).
\end{align*}
This shows that $\{\|\nabla u_n\|_2\}$ is  away from zero.  Define  \(\widetilde{K} := \lim\limits_{n\to\infty}\|\nabla u_n\|_2^2 > 0\).   Passing to the limit in \eqref{eq1:Lem2:solution-1} yields
\begin{align*}
&\Bigl[a+b\widetilde{K}+\frac{b^2\mathcal{S}^{3}}{2}+\frac{b\mathcal{S}}{2}\sqrt{b^2\mathcal{S}^4+4(a+b\widetilde{K})\mathcal{S}}\Bigr]\int_{\mathbb{R}^3} \nabla u  \nabla\phi\,\mathrm{d}x \nonumber\\
= & \lambda\int_{\mathbb{R}^3}  u \phi\,\mathrm{d}x \int_{\mathbb{R}^3} |u|^4 u \phi\,\mathrm{d}x + \int_{\mathbb{R}^3} u \phi \log u^2\,\mathrm{d}x.\label{eq2:Lem2:solution-1}
\end{align*}
Hence \(u\) satisfies
\begin{equation*}
-\Bigl[a+b\widetilde{K}+\frac{b^2S^{3}}{2}+\frac{b\mathcal{S}}{2}\sqrt{b^2\mathcal{S}^4+4(a+b\widetilde{K})\mathcal{S}}\Bigr]\Delta u  =  \lambda  u+|u|^4u +u\log u^2.
\end{equation*}
The associated Pohozaev identity reads
\begin{equation}\label{eq3:Lem2:sulotion-11}
\widetilde{P}(u) := \Bigl[a+b\widetilde{K}+\frac{b^2\mathcal{S}^{3}}{2}+\frac{b\mathcal{S}}{2}\sqrt{b^2\mathcal{S}^4+4(a+b\widetilde{K})\mathcal{S}}\Bigr]\|\nabla u\|_2^2 - \|u\|_6^6 - \frac{3 c^2}{2} = 0.
\end{equation}
Setting \(v_n := u_n - u\).  Noting that  \(v_n \rightharpoonup 0\) in \(W_{r}\), \(v_n \to 0\) in \(L^l(\mathbb{R}^3)\) for $l\in(2,6)$ and we see that \eqref{Lem:pro-eq3a11}, \eqref{eq1-lem;duishu} and \eqref{Lem:pro-eq3a1} are still valid.
Using \eqref{eq3:Lem2:sulotion-11} and $P(u_n) \to 0$, we obtain
\begin{equation}\label{eq6:Lem2:sulotion-1}
\|v_n\|_6^6 = \Bigl[a+b\widetilde{K}+\frac{b^2\mathcal{S}^{3}}{2}+\frac{b\mathcal{S}}{2}\sqrt{b^2\mathcal{S}^4+4(a+b\widetilde{K})\mathcal{S}}\Bigr]\|\nabla v_n\|_2^2 + o_n(1),\ \text{as} \  n\to+\infty.
\end{equation}
which implies \eqref{Lem:pro-eq41-1}.
Let \(\ell := \lim\limits_{n\to\infty}\|\nabla v_n\|_2^2\). Using arguments similar to those in Lemma  \ref{prop:PS-critical} either \(\ell = 0 \) or \(\ell \geq M\), where $M$ is defined by \eqref{M}.

\textbf{Step 4: Compactness results for  the Palais-Smale sequences.}
We claim that $l=0$. Indeed, using \eqref{Lem:pro-eq3a11},\eqref{eq1-lem;duishu},\eqref{Lem:pro-eq3a1},\eqref{eq:technical},\eqref{eq6:Lem2:sulotion-1}  and $(iii)$ of Lemma \ref{lem:st},
we have
\begin{align*}
&\ \ \  \ \hat{m} _{{\widetilde{R}}}(c)  =\mathcal{J}(u_n)+o_n(1)\\
&\ge \mathcal{J}(u) +\frac{(b^2 \mathcal{S}^4 + 4a\mathcal{S})^{3/2}}{24} \left[\left(1 + \frac{4b}{b^2 \mathcal{S}^3 + 4a} \|\nabla v_n\|_2^2\right)^{3/2} - 1\right] \\
&\quad
+
 \left(\frac{a}{2} + \frac{b^2 \mathcal{S}^3}{4}\right) \|\nabla v_n\|_2^2+\frac{b}{2}\|\nabla u\|_2^2\|\nabla v_n\|_2^2+ \frac{b}{4} \|\nabla  v_n\|_2^4\\
&\quad - \frac{1}{6} \Bigl[a+b\widetilde{K}+\frac{b^2\mathcal{S}^{3}}{2}+\frac{b\mathcal{S}}{2}\sqrt{b^2\mathcal{S}^4+4(a+b\widetilde{K})\mathcal{S}}\Bigr]\|\nabla v_n\|_2^2+o_n(1)\\
&\ge \mathcal{J}(u) + \Bigl[\frac{a}{3}+\frac{b^2\mathcal{S}^3}{6}+\frac{b}{3}\|\nabla u\|_2^2-\frac{b\mathcal{S}}{12} \sqrt{b^2\mathcal{S}^4+4(a+b{\widetilde{R}}^2)\mathcal{S}}\Bigr]\|\nabla v_n\|_2^2
+\frac{b}{12}\|\nabla v_n\|_2^4+o_n(1)\\
&\ge \mathcal{J}(u) +o_n(1)
\end{align*}
as $n\to+\infty$. By the weak lower semi-continuity and Lemma \ref{lem:mAc=m+}, we have  \( u \in V _{{\widetilde{R}}}(c)\). Consequently, \(\hat{m} _{{\widetilde{R}}}(c)  \leq \mathcal{J}( u)\), which yields  \(\ell = 0\). Then \(u_n \to u\) strongly in \(W_{r}\) and \(L^6(\mathbb{R}^3)\) $ \text{as}\  n\to+\infty$. Hence
 \(u\) is an  local minimizer of \(\mathcal{J}|_{V _{{\widetilde{R}}}(c)}\).  
 Using arguments similar to those in Theorem \ref{thm:main1}, we conclude that \(u > 0\)  is  nonnegative, radially symmetric and non-increasing in $r = |x|$.
 Then the conclusion of Lemma \ref{Lem2:solution} follows by taking $(\hat{\lambda}_{c}, \hat{u}_{c})=(\lambda, u)$. The proof is complete.
\end{proof}
 Define
\begin{equation}\label{c11}
c_{1}:= \left\{\frac{2}{K_q C^{q}_{3,q}} \left[ \left(\frac{a}{3} + \frac{b^2 \mathcal{S}^3}{6} + \frac{b\mathcal{S}}{12}\sqrt{b^2 \mathcal{S}^4 + 4(a + b{\widetilde{R}}^2)\mathcal{S}}\right){{\widetilde{R}}}^{2-q\gamma_q} + \frac{b}{12}{{\widetilde{R}}}^{4-q\gamma_q} \right] \right\}^{\frac{1}{q(1-\gamma_q)}}.
\end{equation}
\begin{lemma}\label{lem2:c0c1}
Assume 
$c\in(0,\min\{\widetilde{c}_{**},c_0,c_{1}\})$, where $\widetilde{c}_{**}$, $c_{0}$ and $c_{1}$ are defined  in Lemma \ref{lem:critical_mass}, Lemma
\ref{lem:c_02} and \eqref{c11}, respectively, and the technical condition \eqref{eq:technical} holds. Then $\hat{u}_{c}$ and $\hat{m} _{{\widetilde{R}}}(c)$ given in Lemma \ref{Lem2:solution} satisfy the following properties.
\begin{itemize}
  \item[(i)]  $\|\nabla \hat{u}_{c}\|_2^2 \to 0$
    and  $\hat{m} _{{\widetilde{R}}}(c) \to 0$ as $c \to 0^+$.
  \item[(ii)]  There exists  $c_{2}>0$
 such that if $0<c\le c_2$,   one has
    \begin{align}\label{mmm}
    \hat{m} _{{\widetilde{R}}}(c)\le &\left[\left(\frac{a}{3} + \frac{b^2 \mathcal{S}^3}{6} + \frac{b\mathcal{S}}{12}\sqrt{b^2 \mathcal{S}^4 + 4(a + b{\widetilde{R}}^2)\mathcal{S}}\right){\widetilde{R}}^{2-q\gamma_q} \right. \nonumber\\
& \ \  \left. + \frac{b}{12}{\widetilde{R}}^{4-q\gamma_q}- \frac{ K_q C^{q}_{3,q}}{2}c^{q(1-\gamma_q)}\right]{\widetilde{R}}^{q\gamma_q}.\end{align}
\end{itemize}
\end{lemma}

\begin{proof}
(i) From \cite{Cazenave1983,dAveniaMontefuscoSquassina2014}, we have $w_{0}\in S(c)$, where
\begin{align}\label{gai-9}
  w_{0}(x):=c\pi^{-\frac{3}{4}}\,e^{-\frac{|x|^{2}}{2}}.
\end{align}
  A direct computation gives
\[
\int_{\mathbb{R}^3} |\nabla w_0|^2\,\mathrm{d}x  = \frac{3c^2}{2},
\]
which implies that $w_0\in V _{{\widetilde{R}}}(c)$ from Lemma \ref{lem:c_02}.
Hence, for   sufficiently small $c$,
by using $(i)$ of Lemma \ref{lem:st}, we have \begin{align*}
\hat{m} _{{\widetilde{R}}}(c)\le \mathcal{J}( w_0)
&\le \left[\frac{a}{2}+\frac{b^2\mathcal{S}^3}{4}+\frac{b\mathcal{S}^{\frac{3}{2}}(b^2\mathcal{S}^4+4a\mathcal{S})}{3}\right]\|\nabla w_0\|_2^2 \nonumber\\
&\quad + \frac{b}{4}\|\nabla w_0\|_2^4 + \frac{1}{2}c^2- \frac{1}{6}\|w_0\|_6^6 - \frac{ 1}{2}\int_{\mathbb{R}^{3}}w_{0}^{2}\log w_{0}^{2}\,\mathrm{d}x \nonumber\\
&\leq \left[\frac{a}{2}+\frac{b^2\mathcal{S}^3}{4}+\frac{b\mathcal{S}^{\frac{3}{2}}(b^2\mathcal{S}^4+4a\mathcal{S})}{3}\right]\frac{3c^2}{2} \nonumber\\
&\quad + \frac{b}{4}\cdot\frac{9c^4}{4} 
+ \frac{1}{2}c^2 - \frac{1}{2}\left[c^2\log(c^2) - \frac{3c^2}{2}\log\pi - \frac{3c^2}{2}\right] \nonumber\\
&=- c^2\log c + c^2\left[\frac{3a}{4} +\frac{3b^2\mathcal{S}^3}{8}+\frac{b\mathcal{S}^{\frac{3}{2}}(b^2\mathcal{S}^4+4a\mathcal{S})}{2}+ \frac{5}{4} + \frac{3}{4}\log\pi\right] + \frac{9bc^4}{16}.\label{eqq1}
\end{align*}
It follows that $\hat{m} _{{\widetilde{R}}}(c)\le 0$ as $c \to 0^+$.

From   $P(\hat{u}_{c})=0$ and
\eqref{eq:technical}, we have
\begin{align*}
\hat{m} _{{\widetilde{R}}} (c)&\ge \Bigl[\frac{a}{3}+\frac{b^2\mathcal{S}^3}{6}-\frac{b\mathcal{S}}{12} \sqrt{b^2\mathcal{S}^4+4(a+b\widetilde{R}^2)\mathcal{S}}\Bigr]\|\nabla \hat{u}_{c}\|_2^2\\
&
\quad+\frac{b}{12}\|\nabla \hat{u}_{c}\|_2^4+\frac{1}{2}c^2-\frac{ K_q C^{q}_{3, q} }{2} \|\nabla \hat{u}_{c}\|_2^{q\gamma_q} c^{ (1 - \gamma_q)q}\\
&\ge  \frac{1}{2}c^2-\frac{ K_q C^{q}_{3, q} }{2} {\widetilde{R}}^{q\gamma_q} c^{ (1 - \gamma_q)q}.
\end{align*}
 From the above, one has  $\hat{m} _{{\widetilde{R}}}(c)\geq 0$ as $c \to 0^+$. Combining this and $\hat{m} _{{\widetilde{R}}}(c)\le 0$ as $c \to 0^+$,  we obtain $\hat{m} _{{\widetilde{R}}}(c)\to 0$ as $c \to 0^+$.

 (ii) From the definition of $c_{1}$,  the right hand side of \eqref{mmm}  is a  positive constant.   
  Then the conclusion follows directly from this and $(i)$.  The
proof is complete.
\end{proof}

\subsubsection{Mountain pass geometry and Palais-Smale sequence construction}

Following the methodology developed in \cite{Jeanjean1997}, we introduce the functional
\[
I_-(u) := I(s_u^- * u).
\]
The differentiability of the functional is established in the following lemma, which is crucial for the variational framework.

\begin{lemma}\label{lem:derivative-formula}
Let $ p= 6$ and $c\in(0,\widetilde{c}_{**})$,  where $\widetilde{c}_{**}$ is defined in Lemma~\ref{lem:critical_mass}. Then for any $u \in S(c)$, we have
$$\mathrm{d}I_-(u)[\psi] = \mathrm{d}I(s_u^- * u)[s_u^- * \psi] \quad \text{for all } \psi \in T_u S(c).$$

\end{lemma}

\begin{proof}
For any $u \in S(c)$  and \(\psi \in T_u S(c)\). Consider the functional \(h(t) := \frac{u + t\psi}{\|u + t\psi\|_2}c \in S(c)\) with \(h(0) = u\) and \(h'(0) = \psi\). Set \(s_t := s_{h(t)}^-\) and \(s_0 := s_u^-\). Then we have 
\[
I_-(h(t)) - I_-(h(0)) = I(s_t * h(t)) - I(s_0 * h(0)).
\]
In view of Lemma \ref{lem:critical-points}, the map $u \mapsto s_u^-$ is of class $C^1$, with $s_0$ serving as a  local  maximizer of  $ I(s*u)$.
Hence, for sufficiently small  $t$, we have
\begin{align*}
&\quad \ \ I\left( s_{t}*h(t) \right) - I\left( s_{0}*h(0)\right) \le I\left(s_{t}*h(t)\right) - I\left({s_t*h(0)}\right) \\
&= \frac{ae^{2s_t}}{2}\left(\|\nabla h(t)\|_2^2-\|\nabla h(0)\|_2^2\right) +\frac{be^{4s_{t}}}{4} \left( \|\nabla h(t)\|_{2}^{4} - \|\nabla h(0)\|_{2}^{4} \right)
\\
&\quad + \left(\frac{1}{2} - \frac{3s_t}{2}\right) \left(\|h(t)\|_2^2 - \|h(0)\|_2^2\right) \\
&\quad - \frac{1}{2}\left(\int_{\mathbb{R}^3} |h(t)|^2 \log|h(t)|^2\,\mathrm{d}x - \int_{\mathbb{R}^3} |h(0)|^2 \log|h(0)|^2\,\mathrm{d}x\right)  - \frac{  e^{6s_t}}{6}\left(\|h(t)\|_6^6-\|h(0)\|_6^6\right)\\
&= ae^{2s_t}\int_{\mathbb{R}^3} t \nabla h(\tau_1 t)   \nabla h'(\tau_1 t) \,\mathrm{d}x+ be^{4s_t}\|\nabla h(\tau_2 t)\|_2^2\int_{\mathbb{R}^3} t\nabla h(\tau_2 t)  \nabla h'(\tau_2 t) \,\mathrm{d}x  \\
&\quad + (1 - 3s_t)\int_{\mathbb{R}^3} h(\tau_3 t) h'(\tau_3 t) t\,\mathrm{d}x  - \int_{\mathbb{R}^3} h(\tau_4 t) h'(\tau_4 t) t \log|h(\tau_4 t)|^2\,\mathrm{d}x \\
&\quad - \int_{\mathbb{R}^3} h(\tau_4 t) h'(\tau_4 t) t\,\mathrm{d}x -   e^{6s_t} \int_{\mathbb{R}^3}|h(\tau_5 t)|^{4} h(\tau_5 t) h'(\tau_5 t) t\,\mathrm{d}x
\end{align*}
for some $\tau_1,\dots,\tau_5\in(0,1)$.
Similarly, one can show that
\begin{align*}
&\quad \ \ I\left( s_{t}*h(t) \right) - I\left( s_{0}*h(0)\right) \geq I\left(s_{0}*h(t)\right) - I\left({s_0*h(0)}\right) \\
&= \frac{ae^{2s_{0}}}{2}\left(\|\nabla h(t)\|_2^2-\|\nabla h(0)\|_2^2\right) +\frac{be^{4s_{0}}}{4} \left( \|\nabla h(t)\|_{2}^{4} - \|\nabla h(0)\|_{2}^{4} \right)
\\
&\quad + \left(\frac{1}{2} - \frac{3s_{0}}{2}\right) \left(\|h(t)\|_2^2 - \|h(0)\|_2^2\right) \\
&\quad - \frac{1}{2}\left(\int_{\mathbb{R}^3} |h(t)|^2 \log|h(t)|^2\,\mathrm{d}x - \int_{\mathbb{R}^3} |h(0)|^2 \log|h(0)|^2\,\mathrm{d}x\right)  - \frac{  e^{6s_{0}}}{6}\left(\|h(t)\|_6^6-\|h(0)\|_6^6\right)\\
&= ae^{2s_{0}}\int_{\mathbb{R}^3} t\nabla h(\tau_6 t)  \nabla h'(\tau_6 t) \,\mathrm{d}x+ be^{4s_{0}}\|\nabla h(\tau_7 t)\|_2^2\int_{\mathbb{R}^3} t\nabla h(\tau_7 t)   \nabla h'(\tau_7 t) \,\mathrm{d}x  \\
&\quad + (1 - 3s_{0})\int_{\mathbb{R}^3} h(\tau_8 t) h'(\tau_8 t) t\,\mathrm{d}x - \int_{\mathbb{R}^3} h(\tau_9 t) h'(\tau_9 t) t \log|h(\tau_9 t)|^2\,\mathrm{d}x \\
&\quad - \int_{\mathbb{R}^3} h(\tau_9 t) h'(\tau_9 t) t\,\mathrm{d}x -   e^{6s_{0}} \int_{\mathbb{R}^3}|h(\tau_{10} t)|^{4} h(\tau_{10} t) h'(\tau_{10} t) t\,\mathrm{d}x
\end{align*}
for some $\tau_6,\dots,\tau_{10}\in(0,1)$.
Combining the above analysis, we obtain
\begin{align*}
&\lim_{t\to 0} \frac{I_-(h(t)) - I_-(h(0))}{t} \\
= &a\int_{\mathbb{R}^3} \nabla ({s_u^-}*u) \cdot \nabla ({s_u^-}*\psi)\,\mathrm{d}x +b\|\nabla ({s_u^-}*u)\|_2^2\int_{\mathbb{R}^3} \nabla ({s_u^-}*u) \cdot \nabla ({s_u^-}*\psi)\,\mathrm{d}x\\
&- \int_{\mathbb{R}^3} ({s_u^-}*u) ({s_u^-}*\psi)  \log({s_u^-}*u)^2\,\mathrm{d}x  -  \int_{\mathbb{R}^3}|({s_u^-}*u)|^{4} ({s_u^-}*u) ({s_u^-}*\psi)\,\mathrm{d}x \\
= &\mathrm{d}I(({s_u^-}*u))[({s_u^-}*\psi)]
\end{align*}
for any $u\in S(c)$ and $\psi\in T_uS(c)$. The proof is complete.
\end{proof}

\begin{lemma}\label{lem:shanlups}
Assume  $ p= 6$ and  $c\in(0,\widetilde{c}_{**})$,  where $\widetilde{c}_{**}$ is defined in Lemma~\ref{lem:critical_mass}. Then there exists a Palais-Smale sequence $\{v_{n} \} \subset \Lambda_{c}^{-}$ for $I(u)$ restricted to $S(c) \cap H_{r}^{1}(\mathbb{R}^{3})$ at level $m^-_{c}$ defined by \eqref{gai-6}.
\end{lemma}

\begin{proof}

Let $X := S(c) \cap H_r^1(\mathbb{R}^3)$ and  $\Gamma := \{\{u\} : u \in X\}$ be the family of all singletons in $X$. Since $\partial\Gamma = \emptyset$, $\Gamma$ is a homotopy stable family of compact subsets of $X$ (without boundary).

By Lemmas~\ref{lem:critical-points} and \ref{lem:m+-}, we have
\begin{align*}
\sigma^- := \inf_{ G \in \Gamma} \max_{u \in G} I_{-}(u) &= \inf_{u \in S(c) \cap H_r^1(\mathbb{R}^3)} I_{-}(u) = \inf_{u \in \Lambda_{r,c}^-} I(u) = \inf_{u \in \Lambda_c^-} I(u).
\end{align*}
There exists  a sequence of sets $G_{n} = \{u_n\} \subset \Gamma$ satisfying
\[
\max_{u \in G_{n}} I_{-}(u) < \sigma^-+ \frac{1}{n}.
\]
Set
\[
\eta: [0,1] \times S(c) \to S(c), \quad \eta(t, u) =       ( t s_{u}^{-})*u.
\]
Then, we have
\[
\widetilde{G}_{n} := \eta(\{1\} \times G_{n}) = \{ s_{u_{n}}^{-}*u_{n}\}
\]
belong to $\Gamma$. Lemma~\ref{lem:critical-points} ensures that $\widetilde{G}_{n} \subset 
\Lambda_{c}^{-}$ for each $n \in \mathbb{N}$. 
Any $v \in \widetilde{G}_{n}$ can be represented as $v = s_{u}^{-}*u$ for some $u \in G_{n}$, and therefore  $I_{-}(v) = I_{-}(u)$. It follows that
\[
\max_{v \in \widetilde{G}_{n}} I_{-}(v) = \max_{u \in G_{n}} I_{-}(u) < \sigma^- + \frac{1}{n},
\]
which implies that $\widetilde{G}_{n}$ constitutes another minimizing sequence of sets for $\sigma^-$.  Thus, $\widetilde{G}_{n}$ is a homotopy stable family of compact subsets of $X$.

Denoting by $\|\cdot\|_{*}$ the dual norm on $(T_{u}S(c))^{*}$. By Lemma \ref{min-max principle} we  obtain a minimizing sequence $\{w_{n}\} \subset X$ for $I_{-}$ at level $\sigma^-$ with the properties
\begin{equation}\label{lem:shanlu-1}
\begin{aligned}
&\text{(i)} \quad \lim_{n \to \infty} I_{-}(w_{n}) = \sigma^-, \\
&\text{(ii)} \quad \lim_{n \to \infty}  \|\mathrm{d}I_{-}(w_{n})\|_{*}= 0, \\
&\text{(iii)} \quad \lim_{n \to \infty} \mathrm{dist}_{W_{r}}(w_{n}, \widetilde{G}_{n}) = 0.
\end{aligned}
\end{equation}

Set $v_{n} := {s_{w_{n}}^{-}} *w_{n}\in \Lambda_{c}^{-}$. 
By Lemma~\ref{lem:tangent}, the map $T_{w_{n}}S(c) \to T_{v_{n}}S(c)$ defined by $\psi \mapsto (- s_{w_{n}}^{-})*\psi$ constitutes an isomorphism. This together with Lemma \ref{lem:derivative-formula} give
\begin{equation}\label{G-1}
\|\mathrm{d}I|_{S(c)}(v_{n})\|_{*} := \sup_{\|\psi\| \le 1, \psi \in  T_{v_{n}}S(c)} |\mathrm{d}I(v_{n})[\psi]| = \sup_{\|\psi\| \le 1, \psi \in T_{v_{n}}S(c)} |\mathrm{d}I_{-}(w_{n})[(-s_{w_{n}}^{-})*\psi]|.
\end{equation}
We claim that the sequence  $\{s_{w_{n}}^{-}\}$ is   uniformly  bounded away from zero.
Indeed, from $I(v_{n}) = I_{-}(w_{n}) \to \sigma^- = m^{-}_c$ and arguments similar to those in Lemma \ref{lem:m+-achieved}, we deduce that $\{v_{n}\}$ is  uniformly  bounded in $W_{r}$. This implies that $\{\widetilde{G}_{n}\}$ is  uniformly bounded in $W_{r}$, and then by virtue of \eqref{lem:shanlu-1}-(iii), we obtain $\sup_{n} \|\nabla w_{n}\|_{2}^{2} < \infty$. Since each $\widetilde{G}_{n}$ is compact, there exists $\tilde{v}_{n} \in \widetilde{G}_{n}$ such that $\mathrm{dist}_{W_{r}}(w_{n}, \widetilde{G}_{n}) = \|\tilde{v}_{n} - w_{n}\|_{W_{r}}\to 0$ as $n\to \infty$.
From  $\tilde{v}_n \in \Lambda_c^-$, we have $\|\nabla \tilde{v}_n\|_2$ is  uniformly  bounded away from zero, which yields  
\[
\|\nabla w_n\|_2^2 \ge \|\nabla \tilde{v}_n\|_2^2 - \|\nabla(\tilde{v}_n - w_n)\|_2^2 \ge C>0.
\]
Then the boundedness of $\{s_{w_{n}}^{-}\}$ follows from \[
e^{2s_{w_{n}}^{-}} = \frac{\|\nabla v_{n}\|_{2}^{2}}{\|\nabla w_{n}\|_{2}^{2}}.
\] 
From this, one has  $\|(-s_{w_{n}}^{-})*\psi\|_{W_{r}} \le C$. Together this and \eqref{G-1}, we have 
\[
\|\mathrm{d}I|_{S(c)}(v_{n})\|_{*} \to 0 \quad \text{as} \ n \to \infty.
\]
Thus $\{v_{n}\} \subset \Lambda_{c}^{-}$ is   a Palais-Smale sequence for $I$ restricted to $S(c)$ at level $\sigma^-=m^{-}_{c}$. The proof is complete.
\end{proof}

\subsubsection{Test function construction and energy level}

Let us recall some well-known facts on the Talenti functions, which will help us to estimate the mountain pass level of $I$. For any $\varepsilon > 0$, define
\[
u_\varepsilon(x) = 3^{\frac{1}{4} }\frac{\varepsilon^{\frac{1}{2} }}{(\varepsilon^2 + |x|^2)^{\frac{1}{2}}}, \quad x \in \mathbb{R}^3.
\]
Then $u_\varepsilon(x)$ is a solution to the critical problem
\[
-\Delta u = u^5, \quad x \in \mathbb{R}^3,
\]
and it satisfies
\[
\|\nabla u_\varepsilon\|_2^2 = \|u_\varepsilon\|_6^6 = \mathcal{S}^{\frac{3}{2}},
\]
where
\begin{align}\label{gai-4}
\mathcal{S} = \inf_{u \in H^1(\mathbb{R}^3) \setminus \{0\}} \frac{\|\nabla u\|_2^2}{\|u\|_6^2} = \frac{\|\nabla u_\varepsilon\|_2^2}{\|u_\varepsilon\|_6^2}.
\end{align}

\begin{lemma}\label{lem:talenti}
Let $\tau \in C_0^\infty(\mathbb{R})$ be a cut-off function such that $0 \le \tau(x) \le 1$ in $\mathbb{R}$, and
\[
\tau(x) =
\begin{cases}
1, & |x| < R, \\
0, & |x| > 2R,
\end{cases}
\]
where $R>0$ is a constant.
Set $U_\varepsilon(x) = \tau(x) u_\varepsilon(x)$. Then
\begin{align*}
\|\nabla U_\varepsilon\|_2^2 &= \mathcal{S}^{3/2} + O(\varepsilon), \\
\|U_\varepsilon\|_6^6 &= \mathcal{S}^{3/2} + O(\varepsilon^3), \\
\|U_\varepsilon\|_q^q &=
\begin{cases}
O(\varepsilon^{q/2}), & 1 \le q < 3, \\
O(\varepsilon^{q/2} |\ln \varepsilon|), & q = 3, \\
O(\varepsilon^{3 - q/2}), & 3 < q < 6,
\end{cases}
\end{align*}
and
\begin{align}\label{lem:talenti-eq1}
\int_{\mathbb{R}^{3}} U_\varepsilon^2 \ln U_\varepsilon^2 \,\mathrm{d}x = C \varepsilon \ln \frac{1}{\varepsilon} + O(\varepsilon),
\end{align}
as $\varepsilon \to 0$.
\end{lemma}

\begin{proof}
We only prove \eqref{lem:talenti-eq1},  and the proof of other results is similar and can be referred
to  \cite{{BrezisNirenberg1983}}. According to the definition of $U_\varepsilon$, we have
\begin{align*}
\int_{\mathbb{R}^{3}} U_\varepsilon^2 \ln U_\varepsilon^2 \,\mathrm{d}x
&= \int_{\mathbb{R}^{3}} \tau^{2}(x) u_\varepsilon^{2}(x) \ln \left( u_\varepsilon^2(x) \tau^{2}(x) \right) \,\mathrm{d}x \\
&= \int_{\mathbb{R}^{3}} \tau^{2}(x) u_\varepsilon^{2}(x) \ln \tau^{2}(x) \,\mathrm{d}x
+ \int_{\mathbb{R}^{3}} \tau^{2}(x) u_\varepsilon^{2}(x) \ln u_\varepsilon^2(x) \,\mathrm{d}x \\
&:= J_1 + J_2.
\end{align*}

Using the properties of the cut-off function $\tau$, we have 
\begin{align*}
|J_1| = \left| \int_{B_{2R}(0) \setminus B_R(0)} \tau^{2}(x) u_\varepsilon^{2}(x) \ln \tau^{2}(x) \,\mathrm{d}x \right| 
\leq C \int_{B_{2R}(0) \setminus B_R(0)} u_\varepsilon^{2}(x) \,\mathrm{d}x
 = O(\varepsilon).
\end{align*}
For $J_2$, we split the integral over different regions
\begin{align*}
J_2 &= \int_{B_R(0)} u_\varepsilon^{2}(x) \ln u_\varepsilon^2(x) \,\mathrm{d}x
+ \int_{B_{2R}(0) \setminus B_R(0)} \tau^{2}(x) u_\varepsilon^{2}(x) \ln u_\varepsilon^2(x) \,\mathrm{d}x =: J_{21} + J_{22}.
\end{align*}
Applying  (iv) of Lemma \ref{lem:st} with  $0<\alpha<1$, we have
\begin{align*}
 J_{22} \leq C \int_{B_{2R}(0) \setminus B_R(0)} u_\varepsilon^{2+\alpha}(x) \,\mathrm{d}x = O(\varepsilon).
\end{align*}
By 
  (iv) of Lemma \ref{lem:st}  with $\alpha>\frac{1}{2}$,  we get
\begin{align*}
 J_{21}
&= C \varepsilon^2 \ln \left( \frac{1}{\varepsilon} \right) \int_{B_{R/\varepsilon}(0)} \frac{1}{1 + |y|^2} \,\mathrm{d}y + C \varepsilon^2 \int_{B_{R/\varepsilon}(0)} \frac{1}{1 + |y|^2} \ln \left( C \frac{1}{1 + |y|^2} \right) \,\mathrm{d}y\\
&\le C \varepsilon^2 \ln \left( \frac{1}{\varepsilon} \right)\int_0^{R/\varepsilon} \frac{r^2}{1 + r^2} \,\mathrm{d}r +C \varepsilon^2\int_0^{R/\varepsilon} \frac{r^2}{(1 + r^2)^{1+\alpha}}\,\mathrm{d}r  \\
&= C\varepsilon \ln \left( \frac{1}{\varepsilon} \right)+ O(\varepsilon).
\end{align*}
Therefore, combining all these estimates yields the desired conclusion, and hence the proof is complete.
\end{proof}

\begin{lemma}\label{A} There holds that 
 $$-\int_{\mathbb{R}^3}
   A(|U_\varepsilon|)\,\mathrm{d}x\ge C\int_{\mathbb{R}^3} U_\varepsilon^2\log U_\varepsilon^2\,\mathrm{d}x+O(\varepsilon), \ \text{as} \ \varepsilon \to 0,$$ where $U_\varepsilon$ is defined in Lemma \ref{lem:talenti}.
\end{lemma}
\begin{proof}
 We decompose $B_{2R}(0)$ into three disjoint subsets 
\begin{align*}
\Omega_1 &:= \left\{x\in B_{2R}(0) : 0\le |U_\varepsilon(x)| < e^{-3}\right\},\\
\Omega_2 &:= \left\{x\in B_{2R}(0) : e^{-3} \le |U_\varepsilon(x)| < 1\right\},\\
\Omega_3 &:= \left\{x\in B_{2R}(0) : |U_\varepsilon(x)| \ge 1\right\}.
\end{align*}
Hence, one has
\begin{equation*}
-\int_{B_{2R}(0)} A(|U_\varepsilon|)\,\mathrm{d}x = -\int_{\Omega_1} A(|U_\varepsilon|)\,\mathrm{d}x - \int_{\Omega_2} A(|U_\varepsilon|)\,\mathrm{d}x - \int_{\Omega_3} A(|U_\varepsilon|)\,\mathrm{d}x.
\end{equation*}

  By the definition of $A(s)$ for $0 \le s < e^{-3}$, we have 
\begin{equation}\label{gai-dui-1}
-A(|U_\varepsilon|) = U_\varepsilon^2 \log U_\varepsilon^2 \quad \text{on } \Omega_1.
\end{equation}
 In $\Omega_2$, both $-A(|U_\varepsilon|)$ and $|U_\varepsilon|^2 \log |U_\varepsilon|^2$ are negative. Define the continuous function
\begin{equation*}
f(s) := \frac{-A(s)}{s^2\log s^2} = \frac{-(3s^2 + 4e^{-3}s - e^{-6})}{s^2\log s^2}, \quad s \in [e^{-3}, 1).
\end{equation*}
Since $\log s^2 < 0$ and $A(s) > 0$ for $s \ge e^{-3}$, we have $f(s) > 0$ and $\lim\limits_{s\to 1^-} f(s) = +\infty$. Consequently, $f$  is bounded from below by some $C > 0$. Therefore,
\begin{equation}\label{gai-dui-2}
-A(|U_\varepsilon|) \ge C\, U_\varepsilon^2 \log U_\varepsilon^2 \quad \text{on } \Omega_2.
\end{equation}
Finally, we address $\Omega_3$.  When $|U_\varepsilon| \ge 1$, we have $U_\varepsilon^2 \log U_\varepsilon^2 \ge 0$ while $-A(|U_\varepsilon|) \le 0$. 
According to  the definition of $A(s)$, for $\beta>0$, there exists $C>0$ such that  $-A(s) \ge -C s^{2+\beta}$.  In view of this and Lemma \ref{lem:talenti}, we have 
\begin{equation}\label{gai-dui-3}
-\int_{\Omega_3} A(|U_\varepsilon|)\,\mathrm{d}x \ge -C\int_{\Omega_3}  |U_\varepsilon|^{2+\beta}\,\mathrm{d}x = O(\varepsilon), \ \text{as} \ \varepsilon \to 0,
\end{equation} 
where $0<\beta<1$.
By
  (iv) of Lemma \ref{lem:st}  with $\alpha=\beta$,   one has
\begin{equation}\label{gai-dui-4}
\int_{\Omega_3} U_\varepsilon^2 \log U_\varepsilon^2\,\mathrm{d}x\le C\int_{\Omega_3 } |U_\varepsilon|^{2+\beta}\,\mathrm{d}x=O(\varepsilon), \ \text{as} \ \varepsilon \to 0.
\end{equation}

 Finally, making use of \eqref{gai-dui-1}-\eqref{gai-dui-4}, we deduce
 \begin{align*}
&-\int_{_{B_{2R}(0)}} A(|U_\varepsilon|)\,\mathrm{d}x\\
\ge &  C\int_{\Omega_1\cup\Omega_2} U_\varepsilon^2 \log U_\varepsilon^2\,\mathrm{d}x -\int_{\Omega_3} A(|U_\varepsilon|)\,\mathrm{d}x\\
=&C\int_{\mathbb{R}^3} U_\varepsilon^2\log U_\varepsilon^2\,\mathrm{d}x
-C\int_{\Omega_3} U_\varepsilon^2\log U_\varepsilon^2\,\mathrm{d}x
-\int_{\Omega_3} A( |U_\varepsilon|)\,\mathrm{d}x\\
\ge&C\int_{\mathbb{R}^3} U_\varepsilon^2\log U_\varepsilon^2\,\mathrm{d}x
-C\int_{\Omega_3} |U_\varepsilon|^{2+\beta}\,\mathrm{d}x\\
\ge&C\int_{\mathbb{R}^3} U_\varepsilon^2\log U_\varepsilon^2\,\mathrm{d}x+O(\varepsilon), \ \text{as} \ \varepsilon \to 0.
\end{align*}
The proof is complete.
\end{proof}

\begin{lemma}\label{lem:cross}
Let $\hat{u} \in S(c)$  obtained in Lemma \ref{Lem2:solution} with $\mathcal{J}(\hat{u}) = \hat{m} _{{\widetilde{R}}}(c)$.   For any $\varepsilon>0$, there exists  $n_{\varepsilon}\in \mathbb{N}$ sufficiently large   such that  
\begin{equation}\label{lem:cross-eq1}
\int_{\mathbb{R}^{3}} \hat{u}(x - n_{\varepsilon}e_{1})U_{\varepsilon}(x)\,\mathrm{d}x = O(\varepsilon), \ \text{as} \ \varepsilon \to 0,
\end{equation}
and
\begin{equation}\label{lem:cross-eq2}
 \int_{\mathbb{R}^{3}} \nabla \hat{u}(x - n_{\varepsilon}e_{1})   \nabla U_{\varepsilon}(x)\,\mathrm{d}x = O(\varepsilon), \ \text{as} \ \varepsilon \to 0,
\end{equation}
where $e_{1}=(1,0,0)$.
\end{lemma}

\begin{proof}

We first estimate the term \eqref{lem:cross-eq1}. Since $\hat{u}$ is positive, radially symmetric and non-increasing, Lemma~\ref{lem:radial-decay} with $t = 2$ yields $\hat{u}(x) \leq C|x|^{-\frac{3}{2}}$. Noting that $\operatorname{supp}(U_{\varepsilon}) \subset B_{2R}(0)$ implies $|x - ne_1| \geq n - 2R$ for all $x \in B_{2R}(0)$. From this and by H\"older's inequality, we obtain
\begin{align*}
\int_{\mathbb{R}^{3}} \hat{u}(x - ne_{1})U_{\varepsilon}(x)\,\mathrm{d}x 
&\leq C\int_{B_{2R}(0)} \frac{U_{\varepsilon}(x)}{|x - ne_{1}|^{\frac{3}{2}}}\,\mathrm{d}x \\
&\leq \frac{C}{(n - 2R)^{\frac{3}{2}}}\int_{B_{2R}(0)} U_{\varepsilon}(x)\,\mathrm{d}x \\
&\leq C(n - 2R)^{-\frac{3}{2}}\|U_{\varepsilon}\|_{2}.
\end{align*}
 Thus \eqref{lem:cross-eq1} follows by taking $n_{\varepsilon}$   sufficiently large  such that $(n_{\varepsilon} - 2R)^{-\frac{3}{2}}\le \|U_{\varepsilon}\|_{2}$.

For the term \eqref{lem:cross-eq2}, integration by parts yields 
\begin{align*}
\left| \int_{\mathbb{R}^{3}} \hat{u}(x - ne_{1}) \Delta U_{\varepsilon}(x) \,\mathrm{d}x \right|
&\leq \int_{B_{2R}(0)} |\hat{u}(x - ne_{1})| |\Delta U_{\varepsilon}(x)|\,\mathrm{d}x \\
&\leq C(n - 2R)^{-\frac{3}{2}} \int_{B_{2R}(0)} |\Delta U_{\varepsilon}(x)| \,\mathrm{d}x \\
&\leq C(n - 2R)^{-\frac{3}{2}} \int_{B_{2R}(0)} |U_{\varepsilon}(x)|^{5} \,\mathrm{d}x, 
\end{align*}
 which implies that  \eqref{lem:cross-eq2}  holds  for  $n_{\varepsilon}$   sufficiently large.
The proof is  complete.
\end{proof}

\begin{lemma}\label{lem:shanlushuiping}
For any $c\in(0,\widetilde{c}_{**})$ and $p=6$, we have that
\begin{equation}\label{con:m}
\begin{aligned}
m^{-}_c
&< \hat{m} _{{\widetilde{R}}}(c)+C(a,b,\mathcal{S}),
\end{aligned}
\end{equation}
where $\widetilde{c}_{**}$, $C(a,b,\mathcal{S})$  and $\hat{m} _{{\widetilde{R}}}(c)$  are defined  in Lemmas~\ref{lem:critical_mass}, \ref{lem:fmax} and \eqref{gai-8} respectively.
\end{lemma}

\begin{proof}

For   $\varepsilon > 0$ small, let $n_{\varepsilon}\in \mathbb{N}$ be as given in Lemma~\ref{lem:cross}, 
and define
$$W_{\varepsilon,t}(x) := \hat{u}(x-n_{\varepsilon}e_1) + tU_\varepsilon,$$
where $\mathcal{J}(\hat{u}) = \hat{m} _{{\widetilde{R}}}(c)$. Let
\begin{align}\label{gai--1}\widetilde{W}_{\varepsilon,t} = s^{\frac{1}{2}}W_{\varepsilon,t}(sx), \quad \text{where} \quad s = \frac{\|W_{\varepsilon,t}\|_2}{c}>0.\end{align}
By  a direct computation, we obtain 
\begin{align*}
&\|\nabla \widetilde{W}_{\varepsilon,t}\|_2^2 = \|\nabla W_{\varepsilon,t}\|_2^2, \\
&\|\widetilde{W}_{\varepsilon,t}\|_2^2 = c^2, \\
&\|\widetilde{W}_{\varepsilon,t}\|_6^6 = \|W_{\varepsilon,t}\|_6^6.
\end{align*}
Therefore, $\widetilde{W}_{\varepsilon,t} \in S(c)$.
According to Lemma \ref{lem:critical-points}, there exists  a unique $s^-_{\varepsilon,t} \in \mathbb{R}$ such that $s^-_{\varepsilon,t}*\widetilde{W}_{\varepsilon,t} \in \Lambda^-_c$.
When $t=0$, we have $\widetilde{W}_{\varepsilon,0} =W_{\varepsilon,0} = \hat{u}(x-n_{\varepsilon}e_1)$ with $P(\hat{u}(x-n_{\varepsilon}e_1))=0$ and $\hat{u}(x-n_{\varepsilon}e_1)\in V _{{\widetilde{R}}}(c)$, which implies  $Q(\hat{u}(x-n_{\varepsilon}e_1))<0$.  Combining  this and Lemma \ref{lem:critical-points},  we get $s^-_{\varepsilon,0} > 0$.
It follows from $Q(s^-_{\varepsilon,t}*\widetilde{W}_{\varepsilon,t})=0$ and  H\"older's inequality  that,  for sufficiently small $\varepsilon > 0$, 
\begin{equation*}
\begin{aligned}
  &e^{6s^-_{\varepsilon,t}}\left(\|\hat{u}(x-n_{\varepsilon}e_1)\|_6^6 + t^6\|U_\varepsilon\|_6^6 \right)\\
\leq& a e^{2s^-_{\varepsilon,t}} \|\nabla \widetilde{W}_{\varepsilon,t}\|_2^2 + b e^{4s^-_{\varepsilon,t}} \|\nabla \widetilde{W}_{\varepsilon,t}\|_2^4 \\
\leq &ae^{2s^-_{\varepsilon,t}}\left(\|\nabla \hat{u}(x-n_{\varepsilon}e_1)\|_2^2 + 2t\|\nabla \hat{u}(x-n_{\varepsilon}e_1)\|_2\|\nabla U_\varepsilon\|_2 + t^2\|\nabla U_\varepsilon\|_2^2\right) \\
&+be^{4s^-_{\varepsilon,t}}\left(\|\nabla \hat{u}(x-n_{\varepsilon}e_1)\|_2^2 + 2t\|\nabla \hat{u}(x-n_{\varepsilon}e_1)\|_2\|\nabla U_\varepsilon\|_2 + t^2\|\nabla U_\varepsilon\|_2^2\right)^2,
\end{aligned}
\end{equation*}
where we have  used   $(\alpha+\beta)^\gamma \geq \alpha^\gamma + \beta^\gamma$ for $\alpha, \beta,\geq 0$ and $\gamma\geq 1$ in the first inequality.
Thus  $s^-_{\varepsilon,t} \to -\infty$ as $t \to +\infty$ uniformly for sufficiently small $\varepsilon > 0$. Following  the continuity of the map $u \mapsto s_u^-$, there exists a  $t_\varepsilon > 0$ such that $s^-_{\varepsilon,t_\varepsilon} = 0$. Consequently, $\widetilde{W}_{\varepsilon,t_\varepsilon} \in \Lambda^-_c$, and we have
\begin{equation*}
m^-_{c} \leq I(\widetilde{W}_{\varepsilon,t_\varepsilon}) \leq \sup_{t \geq 0} I(\widetilde{W}_{\varepsilon,t})
\end{equation*}
for sufficiently small $\varepsilon > 0$.

To complete the proof, we need to prove 
\begin{equation}\label{eq:key-inequality}
I(\widetilde{W}_{\varepsilon,t}) <  \hat{m} _{{\widetilde{R}}}(c)+C(a,b,\mathcal{S}), \ \text{for\ all } t > 0
\end{equation}
 uniformly for sufficiently small $\varepsilon $.
We  note that  $I(\widetilde{W}_{\varepsilon,t}) \leq \mathcal{J}(\widetilde{W}_{\varepsilon,t}) \to \hat{m} _{{\widetilde{R}}}(c)$ as $t \to 0$, and $I(\widetilde{W}_{\varepsilon,t}) \to -\infty$ as $t \to +\infty$ uniformly for sufficiently small $\varepsilon $. Thus, there exists $t_0 > 0$ such that \eqref{eq:key-inequality} holds for $0 < t < \frac{1}{t_0}$ and $t > t_0$.

For $\frac{1}{t_0} \leq t \leq t_0$, we define
$$\delta := \frac{2t}{c^2}\int_{\mathbb{R}^3} \hat{u}(x-n_{\varepsilon}e_1)U_\varepsilon\,\mathrm{d}x.$$
By Lemma \ref{lem:cross}, we have $\delta = O(\varepsilon)$ as $\varepsilon \to 0$. From \eqref{gai--1},  {as} $\ \varepsilon \to 0$, we have
$$s^2 = \frac{\|W_{\varepsilon,t}\|_2^2}{c^2} = 1 + \delta + O({\varepsilon}),$$
 which implies
\begin{align}\label{1-22}
s^{-2} = (1+\delta)^{-1}+ O({\varepsilon}) = 1 - \delta + O(\varepsilon).
\end{align}
From the above, we have $s^{-2}=1+O({\varepsilon})$ {as} $\ \varepsilon \to 0$. Thus, $s^{-2}\log s=O({\varepsilon})$ {as} $\ \varepsilon \to 0$.

If $x\in   B_{2R}$, then  for  $\varepsilon > 0$   and $n_{\varepsilon}\in \mathbb{N}$ sufficiently large, we have  
$$
tU_\varepsilon \ge  \frac{\hat{u}(x-n_{\varepsilon}e_1)}{2} \ \ \text{for} \ \ 
\frac{1}{t_0} \leq t \leq t_0 ,
$$
   since $\hat{u}(x-n_{\varepsilon}e_1) \leq C|n_{\varepsilon} - 2R|^{-\frac{3}{2}}$ by  Lemma~\ref{lem:radial-decay}.
   Using  $(iii)$ of Lemma \ref{A-B} in  $ B_{2R}$, one has
   \begin{align*}
  & \int_{\mathbb{R}^3}
   A(|\hat{u}(x-n_{\varepsilon}e_1) + tU_\varepsilon|)\,\mathrm{d}x\\
   =&\int_{\mathbb{R}^3\setminus B_{2R}}
   A(|\hat{u}(x-n_{\varepsilon}e_1)  )\,\mathrm{d}x+\int_{  B_{2R}}
   A(|\hat{u}(x-n_{\varepsilon}e_1) + tU_\varepsilon|)\,\mathrm{d}x\nonumber \\ 
   \le&\int_{\mathbb{R}^3\setminus B_{2R}}
   A(|\hat{u}(x-n_{\varepsilon}e_1)  )\,\mathrm{d}x+\int_{B_{2R}}
   A(|\hat{u}(x-n_{\varepsilon}e_1)| )\,dx+Ct\int_{B_{2R}}
   A(|U_\varepsilon|)\,\mathrm{d}x\nonumber \\
    =&\int_{\mathbb{R}^3}
   A(|\hat{u}(x-n_{\varepsilon}e_1)  )\,\mathrm{d}x +Ct\int_{B_{2R}}
   A(|U_\varepsilon|)\,\mathrm{d}x\nonumber 
      \end{align*}
Using  this, Lemma \ref{A} and  the fact that $B$ is increasing, we have
 \begin{align}
&\int_{\mathbb{R}^3} W_{\varepsilon,t}^2\log W_{\varepsilon,t}^2\,\mathrm{d}x\nonumber \\
\ge&
\int_{\mathbb{R}^3}
   B(|\hat{u}(x-n_{\varepsilon}e_1) |)\,\mathrm{d}x
  -\int_{\mathbb{R}^3}
   A(|\hat{u}(x-n_{\varepsilon}e_1) + tU_\varepsilon|)\,\mathrm{d}x\nonumber \\
   \ge&
\int_{\mathbb{R}^3}
   B(|\hat{u}(x-n_{\varepsilon}e_1) |)\,\mathrm{d}x
  -\int_{\mathbb{R}^3}
   A(|\hat{u}(x-n_{\varepsilon}e_1)| )\,dx-Ct\int_{\mathbb{R}^3}
   A(|U_\varepsilon|)\,\mathrm{d}x\nonumber \\
   = &\int_{\mathbb{R}^3} |\hat{u}(x-n_{\varepsilon}e_1)|^2\log|\hat{u}(x-n_{\varepsilon}e_1)|^2\,\mathrm{d}x -Ct\int_{B_{2R}}
   A(|U_\varepsilon|)\,\mathrm{d}x\nonumber\\
   \ge &\int_{\mathbb{R}^3} |\hat{u}(x-n_{\varepsilon}e_1)|^2\log|\hat{u}(x-n_{\varepsilon}e_1)|^2\,\mathrm{d}x+ Ct\int_{\mathbb{R}^3} U_\varepsilon^2\log U_\varepsilon^2\,\mathrm{d}x+O(\varepsilon).
  \label{eqq2}
\end{align}

By using \eqref{Gai-1} with $Y=\|\nabla \hat{u}(x-n_{\varepsilon}e_{1})\|_{2}$, \eqref{1-22}, \eqref{eqq2} and Lemma \ref{lem:talenti}, we have
\begin{align*}
&\ \ \ \ I(\widetilde{W}_{\varepsilon,t})\\
&=\frac{a}{2}\|\nabla W_{\varepsilon,t}\|_{2}^{2}+\frac{b}{4}\|\nabla W_{\varepsilon,t}\|_{2}^{4}+\frac{s^{-2}}{2}\| W_{\varepsilon,t}\|_{2}^{2}- \frac{s^{-2}\log s}{2}\| W_{\varepsilon,t}\|_{2}^{2}\\
&\quad-\frac{s^{-2}}{2}\int_{\mathbb{R}^{3}}W_{\varepsilon,t}^{2}\log W_{\varepsilon,t}^{2}\,\mathrm{d}x-\frac{1}{6}\|W_{\varepsilon,t}\|_{6}^{6}
\\
&\leq\frac{a}{2}\|\nabla \hat{u}(x-n_{\varepsilon}e_{1})\|_{2}^{2}+at\int_{\mathbb{R}^{3}}\nabla \hat{u}(x-n_{\varepsilon}e_{1}) \nabla U_{\varepsilon}\,\mathrm{d}x+\frac{at^{2}}{2}\|\nabla U_{\varepsilon}\|_{2}^{2}\\
&\quad+\frac{b}{4}\|\nabla \hat{u}(x-n_{\varepsilon}e_{1})\|_{2}^{4}+bt^{2}\Big(\int_{\mathbb{R}^{3}}\nabla \hat{u}(x-n_{\varepsilon}e_{1}) \nabla U_{\varepsilon}\,\mathrm{d}x\Big)^{2}+\frac{bt^{4}}{4}\|\nabla U_{\varepsilon}\|_{2}^{4}\\
&\quad+\frac{bt^{2}}{2}\|\nabla \hat{u}(x-n_{\varepsilon}e_{1})\|_{2}^{2}\|\nabla U_{\varepsilon}\|_{2}^{2}+bt\|\nabla \hat{u}(x-n_{\varepsilon}e_{1})\|_{2}^{2}\int_{\mathbb{R}^{3}}\nabla \hat{u}(x-n_{\varepsilon}e_{1})\nabla U_{\varepsilon}\,\mathrm{d}x\\
&\quad+bt^{3}\|\nabla U_{\varepsilon}\|_{2}^{2}\int_{\mathbb{R}^{3}}\nabla \hat{u}(x-n_{\varepsilon}e_{1}) \nabla U_{\varepsilon}\,\mathrm{d}x\\
&\quad+\Big(\frac{1}{2}\int_{\mathbb{R}^{3}}|u(x-n_{\varepsilon}e_{1})|^{2}\,\mathrm{d}x+\frac{t^{2}}{2}\int_{\mathbb{R}^{3}}|U_{\varepsilon}|^{2}\,\mathrm{d}x\\
&\quad+t\int_{\mathbb{R}^{3}}
 \hat{u}(x-n_{\varepsilon}e_{1})U_{\varepsilon}\,\mathrm{d}x\Big)\left(1-\frac{2t}{c^{2}}\int_{\mathbb{R}^{3}}\hat{u}(x-n_{\varepsilon}e_{1})U_{\varepsilon}\,\mathrm{d}x\right)\\
&\quad- \frac{1}{2}\int_{\mathbb{R}^{3}}W_{\varepsilon,t}^{2}\log W_{\varepsilon,t}^{2}\,\mathrm{d}x\left(1-\frac{2t}{c^{2}}\int_{\mathbb{R}^{3}}\hat{u}(x-n_{\varepsilon}e_{1})U_{\varepsilon}\,\mathrm{d}x\right)\\
&\quad-\frac{1}{6}\int_{\mathbb{R}^{3}}|\hat{u}(x-n_{\varepsilon}e_{1}|^{6}\,\mathrm{d}x-\frac{ t^{6}}{6}\int_{\mathbb{R}^{3}}|U_{\varepsilon}|^{6}\,\mathrm{d}x+ O(\varepsilon)\\
&\leq\frac{a}{2}\|\nabla \hat{u}(x-n_{\varepsilon}e_{1})\|_{2}^{2}+at\int_{\mathbb{R}^{3}}\nabla \hat{u}(x-n_{\varepsilon}e_{1}) \nabla U_{\varepsilon}\,\mathrm{d}x+\frac{at^{2}}{2}|\nabla U_{\varepsilon}|_{2}^{2}\\
&\quad+\frac{b}{4}\|\nabla \hat{u}(x-n_{\varepsilon}e_{1})\|_{2}^{4}+bt^{2}\Big(\int_{\mathbb{R}^{3}}\nabla \hat{u}(x-n_{\varepsilon}e_{1}) \nabla U_{\varepsilon}\,\mathrm{d}x\Big)^{2}+\frac{bt^{4}}{4}\|\nabla U_{\varepsilon}\|_{2}^{4}\\
&\quad+\frac{bt^{2}}{2}\|\nabla \hat{u}(x-n_{\varepsilon}e_{1})\|_{2}^{2}\|\nabla U_{\varepsilon}\|_{2}^{2}+bt\|\nabla \hat{u}(x-n_{\varepsilon}e_{1})\|_{2}^{2}\int_{\mathbb{R}^{3}}\nabla \hat{u}(x-n_{\varepsilon}e_{1}) \nabla U_{\varepsilon}\,\mathrm{d}x\\
&\quad+bt^{3}\|\nabla U_{\varepsilon}\|_{2}^{2}\int_{\mathbb{R}^{3}}\nabla \hat{u}(x-n_{\varepsilon}e_{1}) \nabla U_{\varepsilon}\,\mathrm{d}x\\
&\quad+\Big(\frac{1}{2}\int_{\mathbb{R}^{3}}|\hat{u}(x-n_{\varepsilon}e_{1})|^{2}\,\mathrm{d}x+\frac{t^{2}}{2}\int_{\mathbb{R}^{3}}|U_{\varepsilon}|^{2}\,\mathrm{d}x+t\int_{\mathbb{R}^{3}}\hat{u}(x-n_{\varepsilon}e_{1})U_{\varepsilon}\,\mathrm{d}x\Big)\\
&\quad-\frac{2t}{c^{2}}\int_{\mathbb{R}^{3}}\hat{u}(x-n_{\varepsilon}e_{1})U_{\varepsilon}\,\mathrm{d}x\int_{\mathbb{R}^{3}}|W_{\varepsilon,t}|^{2}\,\mathrm{d}x\\
&\quad-\Big(\int_{\mathbb{R}^{3}}|\hat{u}(x-n_{\varepsilon}e_{1})|^{2}\log|\hat{u}(x-n_{\varepsilon}e_{1})|^{2}\,\mathrm{d}x+Ct\int_{\mathbb{R}^{3}}U_{\varepsilon}^{2}\log U_{\varepsilon} ^{2}\,\mathrm{d}x\Big)\frac{1}{2}\\
&\quad+\frac{ t}{c^{2}}\int_{\mathbb{R}^{3}}\hat{u}(x-n_{\varepsilon}e_{1})U_{\varepsilon}\,\mathrm{d}x\int_{\mathbb{R}^{3}}W_{\varepsilon,t}^{2}\log W_{\varepsilon,t}^{2}\,\mathrm{d}x\\
&\quad-\frac{1}{6}\int_{\mathbb{R}^{3}}|\hat{u}(x-n_{\varepsilon}e_{1}|^{6}\,\mathrm{d}x-\frac{ t^{6}}{6}\int_{\mathbb{R}^{3}}|U_{\varepsilon}|^{6}\,\mathrm{d}x+ O(\varepsilon) \\
&\leq I(\hat{u})+\frac{at^{2}}{2}|\nabla U_{\varepsilon}|_{2}^{2}+\frac{bt^{4}}{4}|\nabla U_{\varepsilon}|_{2}^{4}-\frac{ t^{6}}{6}\int_{\mathbb{R}^{3}}|U_{\varepsilon}|^{6}\,\mathrm{d}x
\\
&\quad+\frac{bt^{2}}{2}\|\nabla \hat{u}(x-n_{\varepsilon}e_{1})\|_{2}^{2}\|\nabla U_{\varepsilon}\|_{2}^{2}+ O(\varepsilon)+R_{\varepsilon,t}\\
&\leq I(\hat{u})+\frac{a\mathcal{S}^{\frac{3}{2}}}{2}t^{2}+\frac{b\mathcal{S}^{3}}{4}t^{4}-\frac{ \mathcal{S}^{\frac{3}{2}}}{6}t^{6}+\frac{bt^{2}}{2}\|\nabla \hat{u}(x-n_{\varepsilon}e_{1})\|_{2}^{2}\mathcal{S}^{\frac{3}{2}}+O(\varepsilon)+R_{\varepsilon,t}\\
&= I(\hat{u})+\Big(\frac{a+b\|\nabla \hat{u}(x-n_{\varepsilon}e_{1})\|_{2}^{2}}{2}t^{2}+\frac{b\mathcal{S}^{\frac{3}{2}}}{4}t^{4}-\frac{1 }{6}t^{6}\Big)\mathcal{S}^{\frac{3}{2}}+O(\varepsilon)+R_{\varepsilon,t}\\
&\le I(\hat{u})+ C(a,b,\mathcal{S}) + \frac{\left(b^2 \mathcal{S}^4 + 4a\mathcal{S}\right)^{\frac{3}{2}}}{24} \left[ \left(1 + \frac{4b}{b^2 \mathcal{S}^3 + 4a}  \|\nabla \hat{u}(x-n_{\varepsilon}e_{1})\|_{2}^{2}\right)^{\frac{3}{2}} - 1 \right] \\&\ \ + \frac{b^2 \mathcal{S}^3}{4} \|\nabla \hat{u}(x-n_{\varepsilon}e_{1})\|_{2}^{2}+O(\varepsilon)+R_{\varepsilon,t}\\
&= \mathcal{J}(\hat{u}(x-n_{\varepsilon}e_{1}))+C(a,b,\mathcal{S})+O(\varepsilon)+R_{\varepsilon,t}\\
&=\hat{m} _{{\widetilde{R}}}(c)+C(a,b,\mathcal{S})+O(\varepsilon)+R_{\varepsilon,t},
\end{align*}
where
\begin{align*}
R_{\varepsilon,t}
&:=at\int_{\mathbb{R}^{3}}\nabla \hat{u}(x-n_{\varepsilon}e_{1}) \nabla U_{\varepsilon}\,\mathrm{d}x+bt^{2}\left(\int_{\mathbb{R}^{3}}\nabla \hat{u}(x-n_{\varepsilon}e_{1}) \nabla U_{\varepsilon}\,\mathrm{d}x\right)^{2}\\
&\quad+bt|\nabla \hat{u}(x-n_{\varepsilon}e_{1})|_{2}^{2}\int_{\mathbb{R}^{3}}\nabla \hat{u}(x-n_{\varepsilon}e_{1}) \nabla U_{\varepsilon}\,\mathrm{d}x+bt^{3}|\nabla U_{\varepsilon}|_{2}^{2}\int_{\mathbb{R}^{3}}\nabla \hat{u}(x-n_{\varepsilon}e_{1}) \nabla U_{\varepsilon}\,\mathrm{d}x\\
&\quad+\frac{t^{2}}{2}\int_{\mathbb{R}^{3}}|U_{\varepsilon}|^{2}\,\mathrm{d}x+t\int_{\mathbb{R}^{3}}\hat{u}(x-n_{\varepsilon}e_{1})U_{\varepsilon}\,\mathrm{d}x-\frac{2t}{c^{2}} \int_{\mathbb{R}^{3}}\hat{u}(x-n_{\varepsilon}e_{1})U_{\varepsilon}\,\mathrm{d}x\int_{\mathbb{R}^{3}}|W_{\varepsilon,t}|^{2}\,\mathrm{d}x\\
&\quad-\frac{1}{2}Ct\int_{\mathbb{R}^{3}}U_{\varepsilon}^{2}\log U_{\varepsilon}^{2}\,\mathrm{d}x+\frac{t}{c^{2}} \int_{\mathbb{R}^{3}}\hat{u}(x-n_{\varepsilon}e_{1})U_{\varepsilon}\,\mathrm{d}x\int_{\mathbb{R}^{3}}W_{\varepsilon,t}^{2}\log W_{\varepsilon,t}^{2}\,\mathrm{d}x\\
&=\Big(at+bt^{2}\int_{\mathbb{R}^{3}}\nabla \hat{u}(x-n_{\varepsilon}e_{1}) \nabla U_{\varepsilon}\,\mathrm{d}x
+bt|\nabla \hat{u}(x-n_{\varepsilon}e_{1})|_{2}^{2}\\
&\quad+bt^{3}|\nabla U_{\varepsilon}|_{2}^{2}\Big)\int_{\mathbb{R}^{3}}\nabla \hat{u}(x-n_{\varepsilon}e_{1}) \nabla U_{\varepsilon}\,\mathrm{d}x\\
&\quad+\frac{t^{2}}{2}\int_{\mathbb{R}^{3}}|U_{\varepsilon}|^{2}\,\mathrm{d}x+t\int_{\mathbb{R}^{3}}\hat{u}(x-n_{\varepsilon}e_{1})U_{\varepsilon}\,\mathrm{d}x -\frac{2t}{c^{2}} \int_{\mathbb{R}^{3}}\hat{u}(x-n_{\varepsilon}e_{1})U_{\varepsilon}\,\mathrm{d}x\int_{\mathbb{R}^{3}}|W_{\varepsilon,t}|^{2}\,\mathrm{d}x\\
&\quad-\frac{1}{2}Ct \int_{\mathbb{R}^{3}}U_{\varepsilon}^{2}\log U_{\varepsilon}^{2}\,\mathrm{d}x+\frac{t}{c^{2}} \int_{\mathbb{R}^{3}}\hat{u}(x-n_{\varepsilon}e_{1})U_{\varepsilon}\,\mathrm{d}x\int_{\mathbb{R}^{3}}W_{\varepsilon,t}^{2}\log W_{\varepsilon,t}^{2}\,\mathrm{d}x.
\end{align*}
For $n$ sufficiently large,  by Lemma \ref{lem:talenti} and Lemma \ref{lem:cross}, we see 
\begin{equation*}\label{eq5.17}
R_{\varepsilon,t}
=O(\varepsilon) +C \varepsilon \ln \varepsilon + O(\varepsilon) <0,\ \text{as}  \ \varepsilon \to 0.
\end{equation*}
It follows from that
\begin{equation*}
\begin{aligned}
I(\widetilde{W}_{\varepsilon,t})
< \hat{m} _{{\widetilde{R}}}(c)+C(a,b,\mathcal{S})
\end{aligned}
\end{equation*}
for $t_{0}\leq t\leq\frac{1}{t_{0}}$.

Therefore, for $\varepsilon$ sufficiently  small, we deduce that
\eqref{eq:key-inequality} holds for all $t\geq 0.$ The proof is  complete.
\end{proof}

\subsubsection{Existence of the second solution}

 With the  estimate of the energy of mountain pass level \eqref{con:m} established, we   now complete the construction of the second solution.

\begin{proof}[Proof of Theorem \ref{thm:main3}]
By Lemma \ref{lem:shanlups}, we assume that there exists a Palais-Smale sequence $\{u_n\} \subset \Lambda_c^{-}$ for $I$ restricted to $S(c) \cap H_{r}^{1}(\mathbb{R}^{3})$ at level $m^-_c$ with \eqref{con:m}. 

As in the proof of Lemma~\ref{prop:PS-critical}, we obtain
\[
u_n \rightharpoonup u \quad \text{weakly in } W_r, \qquad
u_n \to u \quad \text{strongly in } L^2(\mathbb{R}^3), \qquad \text{and}\ \ \ \ u_n \to u \quad \text{a.e. in } \mathbb{R}^3.
\]
Consequently,
 $u_n \to u$ strongly in $L^l(\mathbb{R}^3)$ for all $l \in \left(2, 6\right)$ by Lemma \ref{lem:G-N}. 
Set  \(v_n := u_n - u\) and $\ell := \lim\limits_{n\to\infty}\|\nabla v_n\|_2^2$. 
  We obtain that \eqref{eq:lagrange2}--\eqref{Lem:pro-eq41-1} remain valid, 
And we  still conclude that  either $\ell = 0$ or
$
\ell \geq M,$ where $M$ is defined by \eqref{M}.

 We claim that \(\ell = 0\).
  Indeed, if
  $\ell\geq M,$
there are two cases to consider.

\textbf{Case 1).} $\|\nabla u\|_2^2 < {\widetilde{R}}^2$. 

Then it follows from  \eqref{Lem:pro-eq3a11}-\eqref{Lem:pro-eq41}
that
\begin{align*}
m^-_c&=I(u_n)+o_n(1)\\
&= I(u) + \frac{a+b\|\nabla u\|_2^2}{2}\|\nabla v_n\|_2^2 + \frac{b}{4}\|\nabla v_n\|_2^4 - \frac{1}{6}\|v_n\|_6^6-\frac{1}{2}\int_{\mathbb{R}^3} v_n^2 \log v_n^2\,\mathrm{d}x+o_n(1)\\
&= I(u) +\frac{a+b\|\nabla u\|_2^2}{3}\|\nabla v_n\|_2^2 + \frac{b}{12}\|\nabla v_n\|_2^4-\frac{1}{2}\int_{\mathbb{R}^3} v_n^2 \log v_n^2\,\mathrm{d}x+o_n(1)\\
&\ge I(u)+ \frac{a + b\|\nabla u\|_2^2}{3}\ell  + \frac{b}{12}\ell ^2 + o_n(1) \\
&\geq  I(u)+\frac{(a + b\|\nabla u\|_2^2)\mathcal{S}}{6}\left[b\mathcal{S}^2 + \sqrt{b^2\mathcal{S}^4 + 4(a + b\|\nabla u\|_2^2)\mathcal{S}}\right] \\
&\quad + \frac{b\mathcal{S}^2}{48}\left[b\mathcal{S}^2 + \sqrt{b^2\mathcal{S}^4 + 4(a + b\|\nabla u\|_2^2)\mathcal{S}}\right]^2  + o_n(1)  \\
&=  I(u)+\frac{ab\mathcal{S}^3}{4} + \frac{b^3\mathcal{S}^6}{24} + \frac{[b^2\mathcal{S}^4 + 4(a + b\|\nabla u\|_2^2)\mathcal{S}]^{3/2}}{24} + \frac{b^2\mathcal{S}^3}{4}\|\nabla u\|_2^2  + o_n(1) \\
&=  I(u) + \frac{(b^2\mathcal{S}^4 + 4a\mathcal{S})^{3/2}}{24}\left[\left(1 + \frac{4b}{b^2\mathcal{S}^3 + 4a}\|\nabla u\|_2^2\right)^{3/2} - 1\right] +\frac{b^2\mathcal{S}^3}{4}\|\nabla u\|_2^2 \\
&\quad+\frac{ab\mathcal{S}^3}{4} + \frac{b^3\mathcal{S}^6}{24} + \frac{(b^2\mathcal{S}^4 + 4a\mathcal{S})^{3/2}}{24}  + o_n(1) \\
&= \mathcal{J}(u)+C(a,b,\mathcal{S})  + o_n(1)\\
  &\geq  \hat{m} _{{\widetilde{R}}} (c)+C(a,b,\mathcal{S})+ o_n(1),
\end{align*}
which contradicts with \eqref{con:m}.

\textbf{Case 2).} $\|\nabla u\|_2^2 \geq {\widetilde{R}}^2$.

 Then using \eqref{Lem:pro-eq3a11}-\eqref{Lem:pro-eq41},   property (ii) of $B(u)$ and (ii) of in Lemma \ref{lem2:c0c1}, one has
\begin{align*}
m^-_c&=I(u_n)+o_n(1)\\
&= \frac{a}{3}\|\nabla v_n\|_2^2 + \frac{b}{12}\|\nabla v_n\|_2^4 + \frac{b}{6}\|\nabla u\|_2^2 \|\nabla v_n\|_2^2  \\&\quad+\frac{a}{3}\|\nabla u\|_2^2 + \frac{b}{12}|\nabla u|_2^4 +\frac{3}{4}c^2- \frac{1}{2}\int_{\mathbb{R}^{3}}u^{2}\log u^{2}\,\mathrm{d}x-\frac{1}{2}\int_{\mathbb{R}^3} v_n^2 \log v_n^2\,\mathrm{d}x + o_n(1) \\
&= \frac{2a + b\|\nabla u\|_2^2}{6} \ell + \frac{b}{12}\ell^2 + \frac{a}{3}\|\nabla u\|_2^2 + \frac{b}{12}\|\nabla u\|_2^4 +\frac{3}{4}c^2- \frac{1}{2}\int_{\mathbb{R}^{3}}u^{2}\log u^{2}\,\mathrm{d}x \\
&\quad -\frac{1}{2}\int_{\mathbb{R}^3} v_n^2 \log v_n^2\,\mathrm{d}x+o_n(1)  \\
&\ge \frac{(2a + b\|\nabla u\|_2^2)\mathcal{S}}{12}\left[bS^2 + \sqrt{b^2 \mathcal{S}^4 + 4(a + b\|\nabla u\|_2^2)\mathcal{S}}\right]  \\
&\quad + \frac{b\mathcal{S}^2}{48}\left[b\mathcal{S}^2 + \sqrt{b^2 \mathcal{S}^4 + 4(a + b\|\nabla u\|_2^2)\mathcal{S}}\right]^2  \\
&\quad + \frac{a}{3}\|\nabla u\|_2^2 + \frac{b}{12}\|\nabla u\|_2^4 +\frac{3}{4}c^2- \frac{1}{2}\int_{\mathbb{R}^{3}}u^{2}\log u^{2}\,\mathrm{d}x +o_n(1)  \\
&= \frac{ab\mathcal{S}^3}{4} + \frac{b^3 \mathcal{S}^6}{24} + \frac{b^2 \mathcal{S}^4 + 2(2a + b\|\nabla u\|_2^2)\mathcal{S}}{24}\sqrt{b^2 \mathcal{S}^4 + 4(a + b\|\nabla u\|_2^2)\mathcal{S}}  \\
&\quad + \left(\frac{a}{3} + \frac{b^2 \mathcal{S}^3}{6}\right)\|\nabla u\|_2^2 + \frac{b}{12}\|\nabla u\|_2^4  +\frac{3}{4}c^2- \frac{1}{2}\int_{\mathbb{R}^{3}}u^{2}\log u^{2}\,\mathrm{d}x +o_n(1) \\
&\ge \frac{ab\mathcal{S}^3}{4} + \frac{b^3 \mathcal{S}^6}{24} + \frac{(b^2 \mathcal{S}^4 + 4a\mathcal{S})^{3/2}}{24} \\
&\quad + \left(\frac{a}{3} + \frac{b^2 \mathcal{S}^3}{6} + \frac{b\mathcal{S}}{12}\sqrt{b^2 \mathcal{S}^4 + 4(a + b\|\nabla u\|_2^2)\mathcal{S}}\right)\|\nabla u\|_2^2 + \frac{b}{12}\|\nabla u\|_2^4 \\
&\quad - \frac{ K_q C^{q}_{3,q}}{2}\|\nabla u\|_2^{q\gamma_q}c^{q(1-\gamma_q)} + o_n(1) \\
&\ge C(a,b,\mathcal{S}) + \left[\left(\frac{a}{3} + \frac{b^2 \mathcal{S}^3}{6} + \frac{b\mathcal{S}}{12}\sqrt{b^2 \mathcal{S}^4 + 4(a + b{\widetilde{R}}^2)\mathcal{S}}\right){\widetilde{R}}^{2-q\gamma_q} + \frac{b}{12}{\widetilde{R}}^{4-q\gamma_q}\right. \\
&\qquad \left. - \frac{ K_q C^{q}_{3,q}}{2}c^{q(1-\gamma_q)}\right]\|\nabla u\|_2^{q\gamma_q} +o_n(1)\\ 
&\ge C(a,b,\mathcal{S}) + \left[\left(\frac{a}{3} + \frac{b^2 \mathcal{S}^3}{6} + \frac{b\mathcal{S}}{12}\sqrt{b^2 \mathcal{S}^4 + 4(a + b{\widetilde{R}}^2)\mathcal{S}}\right){\widetilde{R}}^{2-q\gamma_q} + \frac{b}{12}{\widetilde{R}}^{4-q\gamma_q}\right. \\
&\qquad \left. - \frac{ K_q C^{q}_{3,q}}{2}c^{q(1-\gamma_q)}\right]{\widetilde{R}}^{q\gamma_q} +o_n(1)\\
&\ge C(a,b,\mathcal{S})+\hat{m} _{{\widetilde{R}}}(c)+o_n(1),
\end{align*}
which contradicts with \eqref{con:m}.

Consequently, \(\ell = 0\). It follows that \(u_n \to u\) strongly in \(W_{r}\) and  $A(|u_n|) \to A(|u|)$ strongly in $L^1({\mathbb{R}^{3}} )$. Then $u_n \to u$ strongly in $W_{r}$. Finally, using arguments similar to those in Theorem \ref{thm:main1},  we conclude that
there exists a $\lambda\in\mathbb{R} $ such that $(\lambda,u)$ solves  problem \eqref{eq:main-1}-\eqref{eq:main-2}, and $u$ is positive and radially symmetric. Then the conclusion of Theorem \ref{thm:main4} follows by taking $(\lambda_c^-, u_c^-)=(\lambda, u)$. 
The proof is complete.
\end{proof}

\section{Asymptotic analysis}

In this section, we   investigate the limiting behavior of the ground states and the second solutions established in Theorems~\ref{thm:main1}--\ref{thm:main4} as $c \to 0^+$. We first analyze the  boundedness of $m _{R_{*}}(c)$ and $m(c)$.

\begin{lemma}\label{lem:6.2}
\begin{itemize}
\item[(i)]
  Assume  $c \in (0, +\infty)$ if $2\leq p < \frac{14}{3}$, or $c \in (0, \widetilde{c})$ if $p =  \frac{14}{3}$, where $\widetilde{c}$ is defined in Lemma \ref{lem:m(c)-bounded-below}. Then
  \begin{equation*}\label{eq:lem6.2-i}
    m(c)\le
    - c^2\log c + c^2\left(\frac{3a}{4} + \frac{5}{4} + \frac{3}{4}\log\pi\right) + \frac{9bc^4}{16}.
  \end{equation*}
\item[(ii)]
  Assume $\frac{14}{3}<p\leq 6$,   and  $c\in(0,\min\{c_{*},c_0 \})$, where $c_{*}$ and $c_0$   are defined  in Lemmas~\ref{lem:critical_mass} and
\ref{lem:c_02}, respectively. Then
  \begin{equation*}\label{eq:lem6.2-ii}
    m _{{R}_{*}}(c)\le
    - c^2\log c + c^2\left(\frac{3a}{4} + \frac{5}{4} + \frac{3}{4}\log\pi\right) + \frac{9bc^4}{16}.
  \end{equation*}
\end{itemize}
\end{lemma}

\begin{proof}
(i)   From $w_{0}\in S(c)$, a direct computation gives
\begin{align*}
m(c) \le I(w_{0})&\leq \frac{a}{2}\cdot\frac{3c^2}{2} + \frac{b}{4}\cdot\frac{9c^4}{4} + \frac{1}{2}c^2 - \frac{1}{2}\left[c^2\log(c^2) - \frac{3c^2}{2}\log\pi - \frac{3c^2}{2}\right] \\
&=- c^2\log c + c^2\left[\frac{3a}{4} + \frac{5}{4} + \frac{3}{4}\log\pi\right] + \frac{9bc^4}{16},
\end{align*}
where  $w_{0}$ defined by \eqref{gai-9}.

(ii) 
According to Lemma \ref{lem:c_02}, we have $w_0\in V _{R_{*}}(c)$. Similar to the proof of $(i)$, we obtain the desired conclusion.
 The
proof is complete.
\end{proof}

\textbf{Proof of Theorem \ref{th:to}}
(i) If $2< p\leq 6$, it follows from Theorems \ref{thm:main1}-\ref{thm:main3} that there exists a $u_c^+ \in S(c)$ such that $I(u_c^+) = m(c)$ if $2< p < \frac{14}{3}$ and $I(u_c^+) = m^+_{c}$ if $\frac{14}{3}\leq p \leq 6$ for   small $c > 0$. 
We conclude from Lemma \ref{lem:6.2} that $m^+_{c} \leq 0$ and $m(c) \leq 0$ as  $c \to 0^+$, which implies that $\{u_c^+\} \subset W_{r}$ is uniformly bounded for $c$ small enough. From this and $Q(u_{c}^{+})=0$, we have
 $$m^+_{c} \ge \frac{a}{3} \|\nabla u_{c}^{+}\|_{2}^{2} + \frac{b}{12} \|\nabla u_{c}^{+}\|_{2}^{4}   + \frac{3}{4} c^{2} -\frac{ K_q C^{q}_{3, q }}{2} \|\nabla u_{c}^{+}\|_2^{q\gamma_q} c^{ (1 - \gamma_q)q}\ge  0,$$
  and $m(c)\geq 0$ as $c \to 0^+$. Hence, $m(c)\to 0$ and $m^+_{c}\to 0$ as $c \to 0^+$. It follows that  $\|\nabla u_c^+\|_2^2 \to 0$.

(ii) If $\frac{14}{3}< p <6$,  it follows from Theorem  \ref{thm:main2}  that there exists a $u_c^- \in \Lambda_{c}^-$ such that $I(u_c^-) = m^-_{c}$. By the Gagliardo-Nirenberg inequality and $\Phi''_{u_c^-}(0)<0$,   we have
\[
2a  \|\nabla u_c^-\|_2^2 + 4b \|\nabla u_c^-\|_2^4<  \gamma_p^2p  \| u_c^-\|_p^p\le \gamma_p^2pC_{3,p}^p\,c^{p(1-\gamma_p)}\|\nabla u_c^-\|_2^{p\gamma_p}, 
\]
which implies 
\[
0\le\frac{2a}{\|\nabla u_c^-\|_2^{p\gamma_p-2}}   + \frac{4b }{\|\nabla u_c^-\|_2^{p\gamma_p-4}} <    \gamma_p^2pC_{3,p}^p\,c^{p(1-\gamma_p)}.
\]
Then, $\|\nabla u_c^-\|_2\to +\infty$ as $c\to0$ since $p\gamma_p>4$.

(iii) If $p = 6$, by Theorem
\ref{thm:main4} and Lemma \ref{lem:shanlushuiping},
 there exists a  $u_c^- \in \Lambda_c^-$ such that
\[
m_{c}^{-}=I(u_c^-) < \hat{m}_{{\widetilde{R}}}(c)+\frac{ab\mathcal{S}^3}{4} + \frac{b^3 \mathcal{S}^6}{24} + \frac{(b^2 \mathcal{S}^4 + 4a\mathcal{S})^{3/2}}{24}.
\]
From $(i)$ of Lemma \ref{lem2:c0c1}, one has
\begin{align}
I(u_c^-) &\leq \frac{ab\mathcal{S}^3}{4} + \frac{b^3 \mathcal{S}^6}{24} + \frac{(b^2 \mathcal{S}^4 + 4a\mathcal{S})^{3/2}}{24}+ o(c),\ \text{as}\  c\to 0.\label{7.1}
\end{align}
By  the Gagliardo-Nirenberg  inequality,   property (ii) of $B(u)$ and $Q(u_{c}^{-})=0$, we have 
\begin{align*}
m_{c}^{-}= I(u_{c}^{-})-\frac{1}{6}Q(u_{c}^{-}) \geq  \frac{1}{3}   a\|\nabla u_{c}^{-}\|_2^2 +\frac{1}{12}   b\|\nabla u_{c}^{-}\|_2^4   -\frac{ K_q C^{q}_{3, q}}{2} \|\nabla u_{c}^{-}\|_2^{q\gamma_q} c^{ (1 - \gamma_q)q}.
\end{align*}
This together with \eqref{7.1} implies  that $\{u_c^-\} \subset W_{r}$ is  uniformly bounded  for $c$ small enough. 
Therefore, up to a subsequence, $\{u_c^-\}$ is weakly convergent in $W_r$, strongly convergent in $L^l(\mathbb{R}^3)$ for $l\in[2,6)$, and convergent almost everywhere in $\mathbb{R}^3$ as $c\to 0$.

Let \(\ell_{0} := \lim\limits_{c\to 0^+}\|\nabla u_c^-\|_2^2\). 
From $Q(u_c^- )=0$ and  the Sobolev inequality, we have 
\begin{equation}\label{Lem:pro-eq3---1}
 a\|\nabla u_c^-\|_2^2+b\|\nabla u_c^-\|_2^4+o(c) = \|u_c^-\|_6^6\le \mathcal{S}^{-3}\|\nabla u_c^-\|_2^6, \ \text{as}\ c\to 0.
\end{equation}
Then we conclude   either \(\ell_{0} = 0\) or
$\ell_{0}\geq M_{0}:= \frac{\mathcal{S}}{2}(b\mathcal{S}^2 + \sqrt{  b^2\mathcal{S}^4+4 a\mathcal{S})}.$ 
From $u_c^- \in \Lambda_c^-$, we have
\[
2a  \|\nabla u_c^-\|_2^2 + 4b \|\nabla u_c^-\|_2^4<    \| u_c^-\|_6^6\le \mathcal{S}^{-3}\|\nabla u_c^-\|_2^6, 
\]
which implies \[
 \| \nabla u_c^-\|_2^4>2a\mathcal{S}^{3}.
\]
Hence, $\ell_{0} \geq M_{0}$. By  the Gagliardo-Nirenberg inequality and property (ii) of $B(u)$,  as $c\to 0$, one has
 \begin{align*}
 \int_{\mathbb{R}^3}  (u_c^{-})^2 \log (u_c^{-})^{2}\,\mathrm{d}x
 \le  \frac{ K_q C^{p}_{3,q}}{2}\|\nabla u_c^{-}\|_2^{q\gamma_q}c^{q(1-\gamma_q)}\to 0.\end{align*}
From this,  $\ell_{0} \geq M_{0}$ and  $Q(u_c^-) = 0$, we have  
\begin{align}
I(u_c^-) 
&\ge \frac{a}{2} \|\nabla u_c^-\|_2^2 + \frac{b}{4}  \|\nabla {u_c^-}\|_2^4 -\frac{1}{6}\| u_c^-\|_6^6+  o(c)\nonumber \\
&\geq \frac{1}{3} a\|\nabla u_c^-\|_2^2 + \frac{1}{12}  b\|\nabla {u_c^-}\|_2^4+ o(c)\nonumber\\
&\ge \frac{a \mathcal{S}}{6}\left[b\mathcal{S}^2 + \sqrt{b^2\mathcal{S}^4 + 4a\mathcal{S}}\right]  + \frac{b\mathcal{S}^2}{48}\left[b\mathcal{S}^2 + \sqrt{b^2\mathcal{S}^4 + 4a \mathcal{S}}\right]^2  + o(c) \nonumber \\
&=\frac{ab\mathcal{S}^3}{4} + \frac{b^3\mathcal{S}^6}{24} + \frac{a\mathcal{S}}{6}\sqrt{b^2\mathcal{S}^4 + 4a \mathcal{S}} +\frac{b^{2}\mathcal{S}^{4}}{24}\sqrt{b^2\mathcal{S}^4 + 4a\mathcal{S}}+ o(c) \nonumber\\
&=\frac{ab\mathcal{S}^3}{4} + \frac{b^3\mathcal{S}^6}{24} + \frac{(b^2\mathcal{S}^4 + 4a\mathcal{S})^{3/2}}{24}+o(c).\label{eqq3}
\end{align}
Combining \eqref{7.1} and \eqref{eqq3}, we get
\begin{equation}\label{1}
I(u_c^-) \to \frac{ab\mathcal{S}^3}{4} + \frac{b^3 \mathcal{S}^6}{24} + \frac{(b^2 \mathcal{S}^4 + 4a\mathcal{S})^{3/2}}{24}=C(a,b,\mathcal{S}) \quad \text{as } c \to 0^+,
\end{equation} i.e., $$I(u_c^-)\to \frac{1}{3} a M_{0}+ \frac{1}{12}  bM_{0}^2, \ \ \text{as} \ \ c\to 0.$$
From this and \eqref{Lem:pro-eq3---1}, we deduce 
\begin{align*}
\|\nabla u_c^-\|_2^2 &\to  M_{0}, \\
\| u_c^-\|_6^6 &\to \mathcal{S}^{-3}M_{0}^{3},
\end{align*}
and
\[m^-_{c}\to C(a,b,\mathcal{S}),
 \\
 \quad \text{as } c \to 0.
\]
Hence,  $\{u_c^-\}$ is a minimizing sequence of the  minimizing problem \eqref{gai-4}.
 The
proof is complete.

\end{document}